\documentclass{amsart}
\usepackage{amsmath, amssymb, amscd, amsthm, amsfonts, amstext, amsbsy,  mathrsfs, hyperref,  upgreek, mathtools, stmaryrd, enumitem, tensor, comment, lineno, tikz-cd}
\hypersetup{colorlinks=true}
\usepackage[textsize=footnotesize, color=blue!40
, bordercolor=blue!80]{todonotes}
\numberwithin{equation}{section}

\renewcommand{\subset}{\subseteq}

\newcommand{\B}{\mathcal{B}}
\newcommand{\C}{\mathcal{C}}
\newcommand{\D}{\mathcal{D}}

\newcommand{\F}{\mathcal{F}}

\newcommand{\I}{\mathcal{I}}
\newcommand{\J}{\mathcal{J}}

\newcommand{\N}{\mathcal{N}}

\newcommand{\U}{\mathcal{U}}

\newcommand{\pre}[2]{\tensor[^{#1}]{#2}{}}
\newcommand{\seq}[2]{\langle #1 \mid #2 \rangle}

\newcommand{\Nbhd}{\boldsymbol{N}}

\newcommand{\markdef}[1]{\textbf{#1}}

\newcommand{\Bor}{\k^+\text{-}\mathbf{Bor}}
\newcommand{\Sii}[2]{\k^+\text{-}\boldsymbol{\Sigma}^{ #1 }_{ #2}}
\newcommand{\Pii}[2]{\k^+\text{-}\boldsymbol{\Pi}^{ #1 }_{ #2 }}
\newcommand{\Dee}[2]{\k^+\text{-}\boldsymbol{\Delta}^{ #1 }_{ #2}}

\newcommand{\Borid}{\I\text{-}\mathbf{Bor}}
\newcommand{\Siid}[2]{\I\text{-}\boldsymbol{\Sigma}^{ #1 }_{ #2}}
\newcommand{\Piid}[2]{\I\text{-}\boldsymbol{\Pi}^{ #1 }_{ #2 }}
\newcommand{\Deid}[2]{\I\text{-}\boldsymbol{\Delta}^{ #1 }_{ #2}}

\renewcommand{\a}{\alpha}
\renewcommand{\b}{\beta}
\newcommand{\g}{\gamma}

\renewcommand{\d}{\delta}
\renewcommand{\k}{\kappa}

\renewcommand{\phi}{\varphi}
\renewcommand{\o}{\omega}

\newcommand{\id}{\operatorname{id}}
\newcommand{\cof}{\operatorname{cof}}

\newcommand{\proj}{\operatorname{p}}
\newcommand{\pow}{\mathscr{P}}

\newcommand{\ZFC}{{\sf ZFC}}

\newcommand{\Ord}{\mathrm{Ord}}

\DeclareMathOperator{\ran}{ran}
\DeclareMathOperator{\dom}{dom}

\DeclareMathOperator{\Fn}{Fn}
\DeclareMathOperator{\CUB}{CUB}
\newcommand{\clopen}[1]{\Nbhd_{#1}}

\newtheorem{theorem}{Theorem}[section]
\newtheorem{lemma}[theorem]{Lemma}
\newtheorem{corollary}[theorem]{Corollary}
\newtheorem{proposition}[theorem]{Proposition}

\theoremstyle{definition}
\newtheorem{claim}{Claim}[theorem]
\newtheorem*{claim*}{Claim}
\newtheorem{definition}[theorem]{Definition}
\newtheorem{fact}[theorem]{Fact}
\newtheorem{example}[theorem]{Example}
\newtheorem{question}[theorem]{Question}
\newtheorem*{question*}{Question}
\theoremstyle{remark}
\newtheorem{remark}[theorem]{Remark}

\newenvironment{enumerate-(1)}{\begin{enumerate}[label={\upshape (\arabic*)}, leftmargin=2pc]}{\end{enumerate}}

\title{On Borel sets in ideal topologies}

\author{Miguel Moreno}
\address{Department of Mathematics and Statistics, University of Helsinki, Pietari Kalmin katu 5, Helsinki 00560, Finland and Mathematics Research Centre, Tampere University, Tietotalo, Korkeakoulunkatu 1, Tampere 33720, Finland. }
\urladdr{http://miguelmath.com}

\author{Beatrice Pitton}
\address[Beatrice Pitton]{
Université de Lausanne, Quartier UNIL-Chamberonne, Bâtiment Anthropole,
1015 Lausanne, Switzerland and Dipartimento di matematica Giuseppe Peano, Universit\`a di Torino, Via Carlo Alberto 10, 10121 Torino, Italy.}
\email{beatrice.pitton@unil.ch}


\begin{document}
\begin{abstract}

We study the Borel and analytic subsets of the spaces \({}^{\kappa}\kappa\) and \({}^{\kappa}2\) endowed with ideal topologies, where \(\kappa\) is a regular uncountable cardinal,  thereby addressing some open problems of the literature. 
We provide a systematic analysis of the Borel hierarchy for an arbitrary ideal topology. In particular, we formulate a sufficient condition ensuring that
the hierarchy does not collapse, demonstrate that every Borel set in such a topology is analytic, and establish the existence of a set that is not Borel. We show that, when the underlying ideal contains an unbounded subset, the collection of analytic sets coincides with the full power set of the ambient space. 
Finally, we prove that the Approximation Lemma holds in the setting of ideal topologies.

\end{abstract}

\maketitle


\section{Introduction}

Where classical descriptive set theory examines the Borel and analytic hierarchies on Polish spaces like the Cantor space and Baire space, generalized descriptive set theory investigates analogous hierarchies on generalized versions of these spaces, obtained by replacing \(\omega\) with an uncountable cardinal \(\kappa\) (or with its cofinality). 
The study of generalized descriptive set theory began with the foundational works of Vaught \cite{Vau75} and Väänänen \cite{MR1110032}, and the literature on \(\pre{\kappa}{2}\) and \(\pre{\kappa}{\kappa}\) for uncountable cardinals \(\kappa\) has since become extensive. 
In most works in this area, these spaces are equipped with the bounded topology under the assumption that \(\kappa^{<\kappa}=\kappa\), yielding interesting results on their \(\kappa^+\)-Borel and \(\kappa\)-analytic subsets (see e.g. \cite{MV, FHK14, HK, LS15, ACMRP}). 
The bounded topology is particularly natural due to its applications in model theory and infinitary logics \cite{Vau75, MV, SV00, She01, SV02, She04, DV11, FHK14, MMR21, Mor23}.

However, meaningful variants of this classical framework have emerged in the literature. One approach is to drop or weaken the cardinal assumption: either working with regular \(\kappa\) while permitting  \(\kappa^{<\kappa}>\kappa\) \cite{HMV}, or adopting the weaker assumption \(2^{<\kappa}=\kappa\) and allowing \(\kappa\) to be singular \cite{AMR, DM, ACMRP,MRP25}. Another direction involves changing the topology itself. While papers on infinitary logics predominantly use the bounded topology, works in general topology \cite{ Kur66, Unl82, Kra01, Ili12, CK13} typically employ either the product or box topology. Recent developments have shown that the product topology can be indispensable in certain contexts in generalized descriptive set theory, as demonstrated in \cite[Sections 12 and 13]{AMR} and \cite[Sections 5 and 6]{DM}.

This brings us to the framework central to this article. In \cite{HKSW22}, Holy, Koelbing, Schlicht, and Wohofsky introduced the \emph{ideal topology} \(\tau_\I\) on \(\pre{\kappa}{2}\) and \(\pre{\kappa}{\kappa}\), where \(\kappa\) is a regular uncountable cardinal and \(\I\) is an ideal extending the bounded ideal (see Definition~\ref{def:ideal_topology}). 
This topology refines the bounded topology and provides a unifying framework that encompasses both the bounded and product topologies as special cases. 
While \cite{HKSW22} studied many properties of ideal topologies and initiated the investigation of \(\I\)-Borel sets, their analysis of the \(\I\)-Borel hierarchy remained incomplete, and they left open fundamental questions about its structure. 
Moreover, the notion of analytic sets was not addressed in the ideal topology context. 
Our work provides a comprehensive treatment of both topics, addressing, or making measurable progress toward resolving, some open problems formulated in \cite{HKSW22}.

We now describe the content of the paper in detail. Section~\ref{sec:ideal} introduces the ideal topology on \(\pre{\kappa}{2}\) and \(\pre{\kappa}{\kappa}\) and establishes its basic properties (Lemma~\ref{lem:basics}). In Section~\ref{sec:continuous}, we characterize continuous functions on these spaces via monotone and domain-increasing functions (Definition~\ref{def_functions}).

Section~\ref{section_hierarchy} provides a comprehensive description of the \(\kappa^+\)-Borel hierarchy on spaces equipped with the ideal topology \(\tau_\I\). 
We provide a rigorous treatment of all relevant definitions and key foundational results about \(\I\)-Borel sets. 
This includes natural features of the \(\I\)-Borel hierarchy, like the property of being increasing (Corollary~\ref{thrm:increasing}), the closure properties of the pointclasses \(\Siid{0}{\a}, \Piid{0}{\a}\) and \(\Deid{0}{\a}\) appearing in it (Proposition~\ref{prop: hierarchy_closure}) and the length of the hierarchy. 

The authors of \cite{HKSW22} identified the following as one of their most important open questions:
\begin{question*}[Question 1, \cite{HKSW22}]
    Do the $\mathcal{I}$-Borel sets form a strict hierarchy of length \(\kappa^+\)?
\end{question*}
We address this question in Theorem~\ref{non-collapse}, establishing that the \(\I\)-Borel hierarchy does not collapse for any ideal \(\I\) that is not tall and extends the bounded ideal.
In this context, we also provide a characterization of the tallness of \(\I\). Specifically, \(\I\) fails to be tall if and only if there exists a homeomorphism between a subspace of \((\pre{\k}{2}, \tau_\I)\) and the space \(\pre{\k}{2}\) endowed with the bounded topology (Theorem~\ref{tall_and_homeo}).

A further question raised in \cite{HKSW22} concerns the potential triviality of the class of \(\I\)-Borel sets. Specifically, the authors noted that, in the case where \(2^{<\k}>\k\) and \(\I\) is stationarily tall, they did not have an example of a set that fails to be \(\I\)-Borel. They asked:
\begin{question*}[Question 2, \cite{HKSW22}]
    Does there always exist a set which is not \(\I\)-Borel (or even does not have the \(\I\)-Baire property)?
\end{question*}
In Theorem~\ref{thrm:withoutBP}, we answer this question in the affirmative by proving that there exists a set that does not have the \(\I\)-Baire property, thus
is not \(\I\)-Borel, including in the previously unresolved case described above.

Section~\ref{section_analytic} introduces and analyzes the class of \(\I\)-analytic sets. Surprisingly, this class behaves quite differently from its classical counterpart. We prove in Theorem~\ref{anal_closure} that, when the ideal \(\I\) contains an unbounded set, the class of \(\I\)-analytic sets coincides with the entire powerset of the space. This dramatic collapse occurs for all reasonable definitions of \(\I\)-analyticity (Corollary~\ref{all_analitic} and Proposition~\ref{prop_wrong_anal}), revealing a fundamental distinction between the bounded topology and ideal topology settings.

Finally, Section~\ref{section_approx} extends the approximation method to all ideal topologies. The Approximation Lemma is one of the principal tools for constructing continuous reductions in generalized descriptive set theory (see \cite{FHK14, HM, FMR21}), and our generalization (Lemma~\ref{approx_lemma}) makes this powerful technique available in the broader context of ideal topologies.

Throughout this paper, we work in \(\ZFC\) and we assume that \(\k\) is an uncountable regular cardinal.

\section{Ideals and ideal topologies}\label{sec:ideal}

An ideal \(\I\) on $\k$ is \markdef{\(\k\)-complete} 
if  for every \(\J \subseteq \I\), \(|\J|<\k\) implies \(\bigcup \J \in \I\), and it is proper if \(\I \neq \pow(\k)\).
Some well known \(\k\)-complete proper ideals are \(b\), the collection of bounded subsets of \(\k\),  and \(\sf{NS_\k}\), the collection of nonstationary subsets of \(\k\).
A set \(\mathcal B \subseteq \I \) is a \markdef{basis for the ideal} \(\I\) if for every \(D \in \I\) there exists \( D' \in \mathcal B\) with \(D \subseteq D'\). For example, \(\{\a \mid \a<\k \}\) is a basis (of size $\k$) for the bounded ideal. Note that there are proper ideals other than the bounded ideal which have a basis of size $\k$, e.g., the ideal generated by the bounded ideal together with a set $A$, where $A$ is an unbounded subset of $\kappa$ with unbounded complement (see Example~\ref{example_cont}). 
Given $\mathcal{I}$ and $\mathcal{J}$ ideals on $\kappa$, we say that $\mathcal{I}$ is \markdef{$\mathcal{J}$-tall} if for all $A\subseteq \kappa$ such that $A\notin\mathcal{J}$, there is $B\in\mathcal{I}\backslash \mathcal{J}$ such that $B\subseteq A$. We say that $\mathcal{I}$ is \markdef{tall} if it is \(b\)-tall, and \markdef{stationary tall} if it is \(\sf{NS_\k}\)-tall.

From this point onward (and
unless stated otherwise), we assume that
\begin{center}
\emph{ $\I$ is a ${}\kappa$-complete proper ideal on $\kappa$ that extends the ideal of bounded sets.}
\end{center}

\begin{definition}\label{approx_property_for_last}
    An ideal $\mathcal{I}$ has the \markdef{approximation property} if there is a strictly increasing continuous{}\footnote{A sequence \(\seq{s_\a}{\a<\k}\) is 
strictly increasing if for every \(\a<\b<\k\), \(s_\a \subsetneq s_\b\), and it is continuous if for every limit ordinal \(\a<\k\), \(s_\a=\bigcup_{\b<\a}s_\b\).} sequence $\seq{ B_\alpha}{\alpha<\kappa}$ in $\mathcal{I}$ 
    such that $B_0=\emptyset$, and for all $B\in \mathcal{I}$ there is $\alpha<\kappa$ such that $B\subsetneq B_\alpha$. We call $\seq{B_\alpha}{\alpha<\kappa}$ an approximation sequence of $\mathcal{I}$.
\end{definition}

\begin{lemma}\label{lem:approximation_property}
    An ideal $\mathcal{I}$ has the approximation property if and only if there is a partition $\{ C_\alpha \mid \alpha<\kappa\}$ of $\kappa$ such that $\{D_\alpha\mid \alpha<\kappa\}$ is a basis for $\mathcal{I}$ and $D_\alpha=\bigcup_{\beta\leq \alpha}C_\beta$.
\end{lemma}
\begin{proof}
    If $\seq{ B_\alpha}{\alpha<\kappa}$ is an approximation sequence of $\mathcal{I}$, then it is a basis for $\mathcal{I}$. Moreover, let $C_0=B_1$ and, for all $0<\alpha<\kappa$, let $C_{\alpha}=B_{\alpha+1}\backslash \bigcup_{\beta< \alpha}C_\beta$. Clearly, $D_\alpha=B_{\alpha+1}$, thus $\{D_\alpha\mid \alpha<\kappa\}$ is a basis for $\mathcal{I}$.

    Vice versa, let $\{ C_\alpha \mid \alpha<\kappa\}$ be a partition of $\kappa$ such that $\{D_\alpha\mid \alpha<\kappa\}$ is a basis for $\mathcal{I}$ and $D_\alpha=\bigcup_{\beta\leq \alpha}C_\beta$. Set $B_\alpha=\bigcup_{\beta<\alpha}D_\alpha$ for every \(\a<\k\). Clearly $\seq{ B_\alpha}{\alpha<\kappa}$ is strictly increasing, and since $\{D_\alpha\mid \alpha<\kappa\}$ is a basis, for all $B\in \mathcal{I}$ there is $\alpha\in \kappa$ such that $B\subseteq D_\alpha\subseteq B_{\alpha+1}$, therefore $\mathcal{I}$ has the approximation property.
\end{proof}

It follows from Lemma~\ref{lem:approximation_property} that $\mathcal{I}$ has the approximation property if and only if it admits a basis of size $\k$.

\medskip
Let \(\nu \in \{2,\k\}\) and consider the set $\pre{\k}{\nu}=\{x \mid x\colon \k \to \nu\}$. 
For any set \(I\subseteq \pow(\k)\), let \(\Fn_I(\pre{\k}{\nu})=\{f \mid f\colon D \to \nu \text{ is a function, and } D \in I\}\).
For \(f\in\Fn_I(\pre{\k}{\nu})\), we define
 \[
 \clopen{f}(\pre{\k}{\nu})=\left\lbrace  x \in \pre{\k}{\nu} \mid  f \subseteq x \right\rbrace. 
 \]
 \begin{definition}\label{def:ideal_topology}
 	Let \(\nu \in \{2,\k\}\). The \markdef{\(\I\)-topology} \(\tau_\I\) on \(\pre{\k}{\nu}\) is the topology generated by the collection
    \[
    \B_\I(\pre{\k}{\nu})=\left\lbrace \clopen{f}(\pre{\k}{\nu})  \mid  f \in \Fn_\I(\pre{\k}{\nu}) \right\rbrace. 
    \]
The family \(\B_\I(\pre{\k}{\nu})\) is the \markdef{canonical basis} for \(\tau_\I\) and its elements are called the (canonical) \markdef{basic \(\I\)-open sets} or \markdef{\(\I\)-cones}.
The elements of \(\tau_\I\) are also called \markdef{\(\I\)-open sets}, and similarly we will address closed (clopen) sets with respect to \(\tau_\I\) as \markdef{\(\I\)-closed} (\markdef{\(\I\)-clopen}) sets.  An \markdef{ideal topology} is an \(\I\)-topology for an ideal \(\I\).
\end{definition}
When the space is clear from the context, we drop it from all the notation above.

Clearly, if \(\I \subseteq \I'\) then \( \tau_\I \subseteq \tau_{\I'}\).  It is immediate that the \(\sf{NS_\k}\)-topology refines the bounded topology. While the \(\I\)-topology coincides with the bounded topology when \(\I=b\), it is strictly finer than \(\tau_{b}\) whenever \(\I \supsetneq b\). Indeed, for any unbounded set \(D\in\I\) and any function \(f:D \to 2\), the set \(\clopen{f}\) is \(\I\)-open and \({b}\)-closed, but not \({b}\)-open. As already mentioned in the introduction, when \(\I=b\) we will refer to $(\pre{\k}{\k}, \tau_b)$ and $(\pre{\k}{2}, \tau_b)$ as the generalized Baire space and the generalized Cantor space, respectively. 
  
The following result is an easy observation about the size of the basic \(\I\)-open sets, which will be useful later on.
\begin{lemma}\label{size_open}
 Let \(\nu \in \{2,\k\}\).  For every \(f \in \Fn_{\I}(\pre{\k}{\nu})\),  \(|\clopen{f}(\pre{\k}{\nu})|>\k\). 
\end{lemma}
\begin{proof}
    Let \(f \in \Fn_\I(\pre{\k}{\nu})\) with \(D=\dom(f)\). Then, \(|\clopen{f}(\pre{\k}{\nu})|=2^{|\kappa\setminus D|}\), so it is sufficient to show that \(\kappa \setminus D\) is unbounded. 
    Towards a contradiction, suppose that 
    \(|\kappa\setminus D|<\kappa\). Then, there exists \(\g<\k\) such that \(\kappa\setminus D \subseteq \g\). Since \(\I\) extends the bounded ideal, \(D \cup \g \in \I \). On the other hand, \(D \cup \g=\k\) hence \(\I\) is not proper, a contradiction.  
\end{proof}
When \(\I=\mathsf{NS}_\k\) and \(\nu=2\), Lemma~\ref{size_open} holds in a stronger form since any \(\I\)-cone with the induced topology is homeomorphic to the space \((\pre{\k}{2}, \tau_{\I}\)) by \cite[Lemma 1.6]{HKSW22}.

\medskip
We recall that for a topological space $(X,\tau)$, the density character is the least possible cardinality of a dense subset of $X$, and the weight is the least possible cardinality of an open basis. Moreover,  $(X,\tau)$ is \(\k\)-compact (or \(\k\)-Lindel\"of) if every \(\tau\)-open covering of $X$ has a subcovering of size smaller than \(\k\).  We say that a subset \(A \subseteq \pre{\k}{\nu}\) is \markdef{\(\I\)-dense} if it is dense with respect to \(\tau_\I\), i.e., if for every \(\I\)-open set \(U \subseteq \pre{\k}{\nu}\), \(A \cap U \neq \emptyset\).  A set \(A \subseteq \pre{\k}{\nu}\) is \markdef{\(\I\)-nowhere dense} if for any \(f \in \Fn_\I(\pre{\k}{\nu})\) there exists \(g \in \Fn_\I(\pre{\k}{\nu})\) such that \(f \subseteq g\) and \(A \cap \clopen{g}=\emptyset\). An \markdef{\(\I\)-meager} set is the union of \(\k\)-many \(\I\)-nowhere dense, and an \markdef{\(\I\)-comeager} set is a set with \(\I\)-meager complement.  Finally, \(A \subseteq \pre{\k}{\nu}\) has the \markdef{\(\I\)-Baire property} if there is an \(\I\)-open set $U$ such that $A \triangle U$ is \(\I\)-meager. We recall that, by the usual argument, every \(\I\)-Borel set has the \(\I\)-Baire property.
We also recall the following result from \cite{HKSW22}, which will be useful later one in Section~\ref{subsec:nontrivial}.
\begin{proposition}\cite[Proposition 1.3]{HKSW22} \label{1.3}The intersection of
\(\k\)-many $\I$-open dense sets is $\I$-dense. Equivalently, for every \(f \in \Fn_\I\), \(\clopen{f}\) is not $\I$-meager. 
\end{proposition}

The next lemma summarizes some properties of the ideal topology on \(\pre{\k}{2}\) and \(\pre{\k}{\k}\). Some of these properties already appeared in \cite{HKSW22}, but we include them for the reader's convenience. 
\begin{lemma}\label{lem:basics}
Consider the spaces \((\pre{\k}{2}, \tau_\I)\) and \((\pre{\k}{\k}, \tau_\I)\). The following properties hold.
\begin{enumerate}[label={\upshape (\arabic*)}, leftmargin=2pc]
    \item\label{T2} The topology \(\tau_\I\) is perfect, regular Hausdorff, and zero-dimensional {}\footnote{A topological space is perfect if it has no isolated points; it is zero-dimensional if it has a basis of clopen sets.}.
    \item\label{add}  The topology \(\tau_\I\) is closed under intersections of length at most \(\a \) (for some ordinal $\a$) if and only if \(\a < \kappa\). Therefore, the collection of all \(\I\)-clopen subsets of \(\pre{\k}{\nu}\) is a \(\k\)-algebra{}\footnote{Let \(\g\) be an ordinal. A \(\g\)-algebra on a set $X$ is a family of subsets of $X$ closed under the operations of complementation and well-ordered unions of length at most \( \g\).}. 
     \item\label{basis} For all \(x  \in \pre{\k}{\nu}\), any \(\I\)-open neighborhood basis of \(x\) has size at least \( \k\). Moreover, each point has an \(\I\)-open neighborhood basis of size \(\k\) if and only if \(\I\) has a basis of size \(\k\).
\end{enumerate}
Moreover, if \(\I\) contains an unbounded subset of \(\k\):
\begin{enumerate}[label={\upshape (\arabic*)}, leftmargin=2pc, resume]
    \item\label{compact} The topology \(\tau_\I\) is not compact nor \(\k\)-compact.
    \item\label{weight} The topology \(\tau_\I\) has weight \(2^\kappa\), and \(|\tau_\I|=2^{2^\k}\).
    \item\label{density} The topology \(\tau_\I\) has density character \(2^\k\).
\end{enumerate}
\end{lemma} 
\begin{proof}
Let \(\nu \in \{2,\k\}\).

\ref{T2} It is easy to check that the sets in \(\B_\I\) are clopen, hence the topology \(\tau_\I\) is zero-dimensional and regular. Moreover, it is perfect and Hausdorff, because \(b \subseteq \I \) and \( \tau_b\) is perfect and Hausdorff. 

\ref{add} Suppose \(\a <\k\), let \(\seq{ U_\b}{\b<\a} \) be a sequence of \(\I\)-open sets, and let \(V=\bigcap_{\b<\a}U_\b\). To show that \(V \) is \(\I\)-open, for every \(x\in V\) we construct  \(f \in \Fn_\I \) such that \(x \in \clopen{f} \subseteq V\). For every \(\b<\a\), let \(D_\b \in \I\) be such that \(\clopen{x\restriction D_\b} \subseteq U_\b\). Then \(D=\bigcup_{\b<\a}D_\b \in \I\) by \(\k\)-completeness of \(\I\), and \(\clopen{x\restriction D} \subseteq V\). 

For all $A\subseteq\kappa$, let $\bar 0_A$ be the function constant to $0$ with domain $A$. Suppose \(\a \geq \k\), and let \(\seq{ D_\b}{\b<\a}\) be an increasing sequence of elements of \(\I\) such that \(\bigcup_{\b<\a}D_\b =\k \notin \I\). Set \(V_\b=\clopen{\bar 0_\k \restriction D_\b}\). Then, \(\bigcap_{\b<\a}V_b=\{\bar 0_{\k}\}\) is an \(\I\)-closed set and it is not \(\I\)-open.

\ref{basis} Let \(x \in \pre{\k}{\nu}\) and \(\U_\I\) be an arbitrary \(\I\)-open neighborhood basis of \(x\). We want to construct a map \(\U_\I \to \Fn_\I\) that sends each \(U \in \U_\I\) to \(f_U \subseteq x\) such that \(\clopen{f_U} \subseteq U \), yielding the neighborhood basis \(\U'=   \left\lbrace \clopen{f_U} \mid U \in \U_\I \right\rbrace \). To do that, for every \(U \in \U_\I\) choose \(D_U \in \I\) such that \(\clopen{x\restriction D_U } \subseteq U \) and consider \(\U'=\left\lbrace \clopen{x\restriction D_U } \mid U \in \U_\I \right\rbrace \). 
Notice that the canonical map \(\U_\I \to \U'\) is a surjection, hence for the first part it is enough to show that \(|\U'|\ge \k\). 
Suppose, towards a contradiction, that \( |\U'|< \k\). Since \(\I\) is \(\k\)-complete, we have \(\bar{D} = \bigcup_{U \in \U_\I} D_U \in \I\). Moreover, for every \(U \in \U_\I\), 
\(\Nbhd_{x\restriction \bar{D}} \subseteq \Nbhd_{x\restriction D_U} \subseteq U\). 
Consequently,  \(\{\Nbhd_{x\restriction \bar{D}}\}\) constitutes an open neighborhood base at \(x\). This contradicts \ref{T2}. 

Next, we take care of the second part the statement. For the implication from left to right, we show that \(\{D_U \mid U \in \U_\I\}\) is a basis for \(\I\), so that \( |\U_\I|= \k\) implies that \(\I\) has a basis of size \(\k\).
Let \(D \in \I\). Clearly, \(\clopen{x \restriction D}\) is an \(\I\)-open neighborhood of $x$, so there exists \(U \in \U_\I\) such that \(x \in U \subseteq \clopen{x \restriction D}\). Thus \(\clopen{x \restriction D_U} \subseteq U \subseteq  \clopen{x \restriction D}\) and \(D \subseteq D_U\).
For the other direction, let \(\B\) be a basis fo \(\I\). It is enough to show that \(\left\lbrace \clopen{x\restriction D}  \mid D \in \B \right\rbrace \) is a neighborhood basis of $x$. Let $V$ be an \(\I\)-open set such that \(x \in V\). There exists \(f \in \Fn_\I\) such that \(x \in \Nbhd_f\subseteq V\). Let \(\dom (f) = D \in \I\), so there exists \(D' \in \B\) such that \(D \subseteq D'\). Therefore \(x \in \clopen{x \restriction D'} \subseteq \clopen{x \restriction D}= \clopen{f}\) and \(\left\lbrace \clopen{x\restriction D}  \mid D \in \B \right\rbrace \) is a neighborhood basis of $x$. 

\smallskip
Suppose now that \(\I\) contains an unbounded subset \(D\) of \(\k\).  Consider the collection of functions \(\F_D=\{f \mid f: D \to \nu\}\), and let \( \C_D= \{\clopen{f} \mid f \in \F_D\}\). Note that \(\C_D\) is an \(\I\)-clopen partition of \(\pre{\k}{\nu}\). Clearly, the sets in \(\C_D\) are \(\I\)-clopen and pairwise disjoint. Moreover, for every \(x \in \pre{\k}{\nu}\), \(x\restriction D=f\) for some \(f \in \F_D\), hence \(x \in \clopen{f}\).

 \ref{compact}  
 Since \( \C_D \) 
 is a \(\I\)-clopen partition of \(\pre{\k}{\nu}\) of size \(2^\k\), no subcover of size \(<\k\) can be extracted.
 
 \ref{weight} Let \(\B\) be an arbitrary basis for \(\tau_\I\). Since $\C_D$ is a partition, the map
\[
\B \to \F_D: B \mapsto \begin{cases}
  \bar{0}_{D} & \text{if } \forall f \in \F_D (B \not\subseteq \clopen{f}) \\
  f & \text{if } B \subseteq \clopen{f}
  \end{cases}
\]
is a well-defined surjection. Since \(|\F_D|=2^\k\), \(\B\ge 2^\k\). On the other hand, \(\{\clopen{f} | \ f \in \Fn_\I\}\) is a basis of size \(2^\k\), hence \(\tau_\I\) has weight \(2^\kappa\).

Clearly, \(|\tau_\I| \leq |\pow(\B)|=2^{2^\kappa}\) as witnessed by the injective map \(\tau_\I \to \pow(\B): U \mapsto \{B \in \B \mid B \subseteq U\}\).
For the other inequality, for all \(X \subseteq \F_D\) define \(U_X=\bigcup_{f \in X}\clopen{f}\). Notice that if $X,Y\subseteq \F_D$, then \(X \neq Y\) implies \(U_X \neq U_Y\). Thus, the map \(\pow(\F_D) \to \tau_\I: X \mapsto U_X\) is injective and we get \(|\tau_\I|\ge2^{2^\k}\). 

    \ref{density} Let \(E \subseteq \pre{\k}{\nu}\) be an \(\I\)-dense set. Since \(\C_D \) is a partition of \(\pre{\k}{\nu}\), the surjection
 \( E\rightarrow \F_D: x \mapsto f \), where  \(f \in \F_D\) is (unique) such that \(x \in \clopen{f}\), witnesses \(|E|\geq 2^\k\). Moreover, for every $f\in \Fn_\I$ define the function $f^0\in \pre{\k}{\nu}$ as follows: $$f^0(\gamma)=\begin{cases} f(\gamma) &\mbox{if }\gamma\in \dom(f)\\
    0 & 
    \text{otherwise.}\end{cases}$$  
 The set \(F^0=\left\lbrace f^0\mid f \in \Fn_\I \right\rbrace \subseteq \pre{\k}{\nu}\) 
  is \(\I\)-dense and has cardinality \(|\Fn_\I|=2^\k\). 
\end{proof}

The next result, though not essential, is of independent interest.
\begin{proposition}
\label{lemma:homeo}
Let \(\nu \in \{2,\k\}\). 
The space \((\pre{\k}{\nu} \times \pre{\k}{\k}, \tau_\I \times \tau_\I)\) is homeomorphic to 
   \((\pre{\k}{\k}, \tau_\I)\).
\end{proposition}
\begin{proof}
     Fix a bijection \(\varphi : \nu \times \k \to \k\), and for every \(\a<\k\), let \((\a)_1<\nu\) and \((\a)_2<\k\) such that \(\phi((\a)_1,(\a)_2)=\a\).
     Let 
    \(
    \Phi:\pre{\k}{\nu} \times \pre{\k}{\k} \to \pre{\k}{\k}
    \) be the function defined by \(  \Phi(x,y)(\a)= \varphi(x(\a),y(\a))\) for every \(\a<\k\).
    We have that \(  \Phi\) is bijective because \(\varphi\) is a bijection. We claim that $\Phi$ is a homeomorphism.
    To see that \(  \Phi\) is continuous, note that for any \(f \in \Fn_\I(\pre{\k}{\k}) \), if \(D=\dom(f)\),
    \begin{align*}
    \Phi^{-1}(\clopen{f})=\{(x,y) \mid \forall \a\in D  \left( x(\a)= (f(\a))_1 \wedge y(\a)=(f(\a))_2\right)\}
    = \clopen{f_1}\times \clopen{f_2},
    \end{align*}
where \(f_1, f_2 \in \Fn_\I\) are defined by \(f_1(\a)=(f(\a))_1\) and \(f_2(\a)=(f(\a))_2\) for every \(\a \in D\).
To see that \(\Phi\) is open, fix some \(f \ \in \Fn_\I(\pre{\k}{\nu})\) and \(g  \in \Fn_\I(\pre{\k}{\k})\), and let \(D_f=\dom(f)\) and \(D_g=\dom(g)\).  Let $x \in \clopen{f}$ and $y \in \clopen{g}$, and define \(h: D_f \cup D_g \to \k \)  by \(h(\a)=\varphi(x(\a),y(\a))\). Clearly $\Phi(x,y)\in \clopen{h}$. Conversely, for all $z\in \clopen{h}$, there are $\zeta_1 \in \pre{\k}{\nu}$ and $\zeta_2 \in \pre{\k}{\k}$ such that $\Phi(\zeta_1,\zeta_2)=z$ and $z\restriction D_f \cup D_g=\Phi(\zeta_1,\zeta_2)\restriction D_f \cup D_g=\Phi(x,y)\restriction D_f \cup D_g$. So for all $\alpha\in D_f\cup D_g$, $\varphi(\zeta_1(\alpha),\zeta_2(\alpha))=\varphi(x(\alpha),y(\alpha))$. Since $\varphi$ is a bijection, $\zeta_1\restriction D_f\cup D_g=x\restriction D_f\cup D_g$ and $\zeta_2\restriction D_f\cup D_g=y\restriction D_f\cup D_g$. Thus $(\zeta_1,\zeta_2)\in \clopen{f}\times\clopen{g}$.
 
\end{proof}

\section{The \(\I\)-continuous functions}
\label{sec:continuous}

Throughout this section, let \(\nu \in \{2,\k\}\). 
We say that a function \(\Phi \colon \pre{\k}{\nu} \to \pre{\k}{\nu}\) is \markdef{\(\I\)-continuous} if it continuous with respect to \(\tau_\I\), i.e., \(\Phi \colon (\pre{\k}{\nu}, \tau_\I) \to (\pre{\k}{\nu}, \tau_\I)\) is continuous. 

Similarly to classical setting (see \cite[Proposition 2.6]{Kech}) and to the generalized setting with respect to the bounded topology \(\tau_{b}\) (see \cite[Definition 3.6]{AMR}), if \(\I\) has the approximation property, then every \(\I\)-continuous function can be represented by a monotone and domain-increasing function  $\varphi \colon \Fn_\I( \pre{\k}{\nu} ) \to \Fn_\I( \pre{\k}{\nu} ) $ (Proposition \ref{prop:approx}). However, when dealing with ideals strictly extending the bounded ideal, the definition is rather different, as limits over direct sets are necessary.

Let \( X \) be a topological space, and let \( (\D,\leq) \) be a directed set, still denoted by \( \D \). A point \( x \in X  \) is a \( \D \)-limit of a \( \D \)-net \( (x_d)_{d \in \D} \) of points from \( X \) if for every open neighborhood \( U \) of \( x \) there is \( d \in \D \) such that \( x_{d'} \in U \) for all \( d' \in \D \) with \( d' \geq d\). When this happens, we write \( x = \lim_{d \in \D} x_d\). 
We will consider \( \D \)-limits where \(\D=(\I, \subseteq)\).

In this section, together with $\I$, also $\J$ is a ${}\kappa$-complete proper ideal on $\kappa$ that extends the ideal of bounded sets.

\begin{definition}\label{def_functions} 
Let  $\varphi \colon \Fn_\J( \pre{\k}{\nu} ) \to \Fn_\I( \pre{\k}{\nu} ) $ be a function. 
We say that $\varphi$ is 
\markdef{monotone} if $ f \subseteq g $ implies $\varphi(f) \subseteq \varphi(g)$ for all $f, g \in \Fn_\J( \pre{\k}{\nu} )$.

We say that $\varphi$ is  \markdef{continuous} if it is monotone and for all $x \in \pre{\k}{\nu} $ and \(D\in \I\) there exists \(E \in \J\) such that \(D \subseteq \dom(\varphi(x \restriction E))\).
\end{definition}

Let  $\varphi \colon \Fn_\J( \pre{\k}{\nu} ) \to \Fn_\I( \pre{\k}{\nu} ) $ be a continuous function. The function induced by $\varphi$ is \(\varphi^* \colon  \pre{\k}{\nu}  \to \pre{\k}{\nu}\) defined for every \(x \in \pre{\k}{\nu}\) by 
\begin{align}\label{eq:cont_star}
\phi^*(x)=\lim_{E \in \J}\phi(x\restriction E).
\end{align}
 Note that since \(\phi\) is continuous, the limit in \eqref{eq:cont_star} exists because for all \(D\in \I\) there exists \(E \in \J\) such that
 \begin{align*}
\phi(x\restriction E')\restriction D=\phi(x\restriction E)\restriction D 
\end{align*}
for all \(E \subseteq E'\) with \(E' \in \J\). Moreover, such limit is unique since \(\pre{\k}{\nu}\) is Hausdorff by Lemma~\ref{lem:basics}\ref{T2}, and therefore \(\phi^*\) is well-defined. 

Notice that if \(\phi\) is continuous, then \(\lim_{E \in \J}\phi(x\restriction E)=\lim_{E \in \B}\phi(x\restriction E)\) for any basis \(\B\) for \(\J\). In particular, when \(\I=\J\) is the bounded ideal, \(\phi^*(x)=\lim_{\a<\k}\phi(x\restriction \a)\) is the usual definition also adopted in \cite[Definition 3.6]{AMR}.

\begin{proposition}\label{prop:cont} 
Let  $\varphi \colon \Fn_\J( \pre{\k}{\nu} ) \to \Fn_\I( \pre{\k}{\nu} )$ be a continuous function. Then, \( \varphi^* \colon (\pre{\k}{\nu}, \tau_\J) \to (\pre{\k}{\nu}, \tau_\I) \) is continuous.
\end{proposition}
\begin{proof}
It is sufficient to check that for every \(f \in \Fn_\I\),
\begin{align*}
(\phi^*)^{-1}(\Nbhd_f)=\bigcup\{\Nbhd_g \mid g \in \Fn_\J, f \subseteq \phi(g)\}.
\end{align*}
Since $\bigcup\{\Nbhd_g \mid  g \in \Fn_\J, f \subseteq \phi(g)\}$ is \(\J\)-open for every \(f \in \Fn_\I\), \(\phi^*\) is continuous. 

First, assume that \(x \in \bigcup\{\Nbhd_g \mid g \in \Fn_\J, f \subseteq \phi(g)\} \) and let \(D=\dom(f)\). Our goal is to show that \(\phi^*(x)\restriction D=f\). Let \(g \in \Fn_\J\) such that \(g \subseteq x \) and \(f \subseteq \phi(g)\). If \(E=\dom(g)\), \(g=x \restriction E\) and \(f = \phi(x \restriction E) \restriction D\). Since \(\phi\) is monotone, for every \(E' \in \J\) such that \(E \subseteq E'\) we have \(\phi(x\restriction E)\restriction D = \phi(x\restriction E') \restriction D\).
By definition of \(\phi^*\), given \(\clopen{\phi^*(x)\restriction D}\)  there exists \(F\in \J\) for every \(F'\in \J\) with \(F \subseteq F'\), \(\phi^*(x)\restriction D = \phi(x \restriction F') \restriction D\). Let \(G= E \cup F \in \J \). Then, for every  \(G' \in \J\) such that \(G \subseteq G'\)
we have 
\[
\phi^*(x) \restriction D= \phi(x \restriction G') \restriction D =f.\] 
Hence \(x \in (\phi^*)^{-1}(\Nbhd_f)\). 

Next, assume that \(x \in (\phi^*)^{-1}(\Nbhd_f)\). Then, \(\phi^*(x)\in \Nbhd_f\) and, if \(D=\dom(f)\), \(\phi^*(x) \restriction D =f\). Since \(\phi\) is continuous,
there exists \(E \in \J \) such that
$D\subseteq \dom(\varphi(x\restriction E))$. Therefore, by definition of \(\phi^*\),
$\varphi(x\restriction E) \restriction D =\phi^*(x)\restriction D$. 
Let \(g=x \restriction E\) to get \(x \in \bigcup\{\Nbhd_g \mid g \in \Fn_\J, f \subseteq \phi(g)\}\).
\end{proof}

One could be tempted to only require that $\varphi$ is monotone and that  \(\bigcup_{D\in \J}\phi(x\restriction D) \in \pre{\k}{\nu}\) for every \(x \in \pre{\k}{\nu}\). However, the function \(\phi^*(x)=\bigcup_{D\in \J}\phi(x\restriction D)\) may fail to be continuous, as the next example shows. 

For all $A\subseteq\kappa$, let $\bar\beta_A$ be the function constant to $\beta$ with domain $A$.

\begin{example}\label{example_cont}
 Let $D\subseteq\k$ be the subset of odd{}\footnote{Recall that every ordinal \(\a\) can be written uniquely as \(\a=\g+n\), with \(n<\o\) and either \(\g=0\) or \(\g\) limit. Accordingly, we say that \(\g\) is even (respectively, odd)
if \(n\) is even (respectively, odd).} ordinals, and let $\I$ the ideal generated by the bounded ideal together with \(D\). Note that \(E \in \I\) if and only if \(E \setminus D \) is bounded. 
 
 Define $\varphi \colon \Fn_\I( \pre{\k}{\nu} ) \to \Fn_\I( \pre{\k}{\nu} )$ as follows. For any \(f \in \Fn_\I( \pre{\k}{\nu} )\), and for every $\alpha\in \dom (f)\setminus D$, set:
 \begin{itemize}
     \item $\varphi(f)(\alpha)=f(\alpha)$,
     \item $\varphi(f)(\alpha+1)=f(\alpha)$. 
 \end{itemize}
 Note that $\dom (\varphi(f))=(\dom (f)\setminus D)\cup \{\alpha+1\mid \alpha\in \dom (f)\setminus D\}$, and since \(\dom(f) \in \I\), we get that \(\dom (f)\setminus D\) is bounded. Therefore, for every $f\in \Fn_\I( \pre{\k}{\nu} )$, $\dom (\varphi(f))$ is bounded.
 
It is easy to verify that \(\phi\) is monotone and that for every \(x \in \pre{\k}{\nu}\), \(\phi^*(x)=\bigcup_{D\in \I}\phi(x\restriction D) \in \pre{\k}{\nu}\). To see this, suppose, towards a contradiction, that there are $E,E'\in \I$ and $\alpha\in \dom(\phi(x\restriction E))\cap \dom(\phi(x\restriction E'))$ such that $\phi(x\restriction E)(\alpha)\neq\phi(x\restriction E')(\alpha)$. Then either $\alpha\notin D$ or $\alpha\in D$. If $\alpha\notin D$, then by the definition of $\phi$, $\phi(x\restriction E)(\alpha)=x(\alpha)=\phi(x\restriction E')(\alpha)$, a contradiction. If $\alpha\in D$, then there is $\beta\notin D$ such that $\alpha=\beta+1$ and $\beta\in \dom(\phi(x\restriction E))\cap \dom(\phi(x\restriction E'))$. By the definition of $\phi$, $\phi(x\restriction E)(\alpha)=x(\beta)=\phi(x\restriction E')(\alpha)$, a contradiction.

However, $\varphi^*$ is not \(\I\)-continuous. 
To see this, let $f=\bar0_D$. We will show that $(\phi^*)^{-1}(\Nbhd_f)$ is not open. 

Notice that $\varphi^*(\bar0_\k)=\bigcup_{E \in \B} \varphi(\bar0_\k \restriction E)=\bar0_\k$. Thus, $\bar0_\k\in (\phi^*)^{-1}(\Nbhd_f)$.
It remains to verify that for all $E\in \I$, $\Nbhd_{\bar0_\k\restriction E}\not\subseteq (\phi^*)^{-1}(\Nbhd_f)$. 

Fix $E\in \I$ and let $x=\bar0_E\cup \bar1_{\k\backslash E}$. Since $\bar0_E=\bar0_\k\restriction E$, $x\in \Nbhd_{\bar0_\k\restriction E}$.  
We have that $\k \setminus (D \cup E)\neq \emptyset$ (otherwise $\k=E\cup D\in \I$ against \(\I\) being proper) so let $\a\in \k \setminus (D \cup E)$. 
Since $\alpha\in \dom (x\restriction E\cup {\alpha+1})\backslash D$, $\alpha+1\in \dom (\varphi(x\restriction E\cup {\alpha+1}))$ and 
\[
\varphi(x\restriction E\cup {\alpha+1})(\alpha+1)=(x\restriction E\cup {\alpha+1})(\alpha)=x(\a)=1.
\]

 Then, 
 \(\varphi^*(x)(\alpha+1)=\phi(x\restriction E\cup {\alpha+1} )(\a+1)=1 \).
Finally, notice that $\alpha+1\in D$, thus $f(\alpha+1)=0\neq 1=\varphi^*(x)(\alpha+1)$. So $x\notin (\phi^*)^{-1}(\Nbhd_f)$.
\end{example}

Proposition~\ref{prop:cont} shows that continuous functions on the space $\Fn_\I( \pre{\k}{\nu} )$ induce \(\I\)-continuous functions. It is not clear whether all the \(\I\)-continuous functions are induced by a continuous functions on the space $\Fn_\I( \pre{\k}{\nu} )$. However, this follows directly from the approximation property.

\begin{proposition}\label{prop:approx}
Assume that \(\I\) has the approximation property. Let  $\Phi \colon \pre{\k}{\nu}  \to \pre{\k}{\nu}$ be an \(\I\)-continuous function. Then, there is a continuous $\varphi \colon \Fn_\I( \pre{\k}{\nu} ) \to \Fn_\I( \pre{\k}{\nu} ) $ such that $\varphi^* = \Phi$. 
\end{proposition}

\begin{proof} 
Let \(\seq{B_\a}{\a<\k}\) be the sequence witnessing that \(\I\) has the approximation property. 
For every \(g \in \Fn_\I\), let \(\a_g=\min \{\a<\k \mid \dom(g)\subseteq B_{\a}\}\), and  define \(X_g=\{f \in \Fn_\I \mid \Phi\left(\Nbhd_g \right) \subseteq \Nbhd_f \text{ and } \dom(f)\subseteq B_{\a_g}\}\). Note that \(X_g\) is 
a lattice and every chain \(\C\) in \(X_g\) has an upper bound (it is sufficient to take \(G=\bigcup_{f \in \C}f\), since $\Phi\left(\Nbhd_g \right)\subseteq\bigcap_{f \in \C}\Nbhd_f=\Nbhd_G$ and \(\dom(G) \subseteq B_{\a_g}\in \I\)).
By Zorn’s lemma, \(X_g\) has at least one maximal element. We now show that this maximal element is unique. Assume that \(f,f'\in X_g\) are both maximal. Since \(\Phi\bigl(\Nbhd_g\bigr) \subseteq \clopen{f} \cap \clopen{f'}\),  \(f\) and \(f'\) must be compatible. Therefore, in order not to violate the maximality of either one, we must have \(f = f'\).
For every \(g \in \Fn_\I\),  let $\varphi(g)$ be the (unique) maximal $f \in X_g$.

Let us show that the function $\varphi$ is monotone. Assume that \(g \subseteq g' \in \Fn_\I\). Note that \(\dom(g) \subseteq B_{\a_g}\) and \(\dom(g) \subseteq \dom(g') \subseteq B_{\a_{g'}}\), so \(\a_g \leq \a_{g'}\).
Moreover, $\Nbhd_{g'}\subseteq \Nbhd_{g}$ and \(\Phi(\Nbhd_{g'})\subseteq \Phi(\Nbhd_{g})\). From the definition of $\phi$ we know that \(\Phi(\Nbhd_{g}) \subseteq \Nbhd_{\phi(g)}\), therefore \(\Phi(\Nbhd_{g'}) \subseteq \Nbhd_{\phi(g)}\). Since \(\dom(\phi(g))\subseteq B_{\a_g}\subseteq B_{\a_{g'}}\), \(\phi(g)\in X_{g'}\). Thus \(\phi(g) \subseteq \phi(g')\), by maximality of \(\phi(g')\).

It remains to prove that \(\phi\) is continuous and that for every $x \in \pre{\k}{\nu}$, $\varphi^*(x)=\Phi(x)$. 
Fix \(x \in \pre{\k}{\nu}\). 
By the \(\I\)-continuity of $\Phi$, for every $D \in \I$ there exists $D'\in \I$ such that $\Phi(\Nbhd_{x \restriction D'}) \subseteq \Nbhd_{\Phi(x)\restriction D}$. 
Let $ E= D \cup D' \in \I$, so $\Nbhd_{x \restriction E}\subseteq \Nbhd_{x \restriction D'}$ and $\Phi(\Nbhd_{x \restriction E})\subseteq \Phi(\Nbhd_{x \restriction D'})$. Then,  $\Phi(\Nbhd_{x \restriction E})\subseteq \Nbhd_{\Phi(x)\restriction D}$, and since \(D \subseteq E \subseteq B_{\a_{x\restriction E }}\), we have
\(\Phi(x)\restriction D\in X_{x\restriction E}\). 
By maximality of \(\varphi(x \restriction E)\), $\Phi(x)\restriction D \subseteq \varphi(x\restriction E)$.
Thus, we have proved that for every $D \in \I$ there exists \( E \in \I\) such that $\Phi(x)\restriction D = \varphi(x\restriction E)\restriction D$. Finally, $\varphi^*(x)=\Phi(x)$ due to the uniqueness of the limit.
\end{proof}

\begin{question}
Does there exist an ideal \(\I\) for which Proposition~\ref{prop:approx} does not hold?
\end{question}

\section{The \(\I\)-Borel sets and their hierarchy}\label{section_hierarchy}
A \markdef{pointclass} \(\boldsymbol\Gamma\) is a class-function assigning to every nonempty topological space $X$ a nonempty family \(\boldsymbol\Gamma(X) \subseteq \pow(X)\). The dual $\check{\boldsymbol\Gamma}$ of \( \boldsymbol{\Gamma}\) is defined by $\check{\boldsymbol\Gamma}(X)= \{X \setminus A \mid A \in \boldsymbol\Gamma(X)\}$. A pointclass \(\boldsymbol{\Gamma}\) is \markdef{boldface} if it is closed under continuous preimages, that is, if $f \colon X \to Y$ is continuous and $B \in \boldsymbol\Gamma(Y)$ implies $f^{-1}(B) \in \boldsymbol\Gamma(X)$.  A pointclass \( \boldsymbol{\Gamma} \) is \markdef{hereditary} if \( \boldsymbol{\Gamma}(Y) = \{ A \cap Y \mid A \in \boldsymbol{\Gamma}(X) \} \) for every \( Y \subseteq X \). Notice that if \( \boldsymbol{\Gamma}\) is a boldface pointclass, then it is hereditary if and only if for every \( Y \subseteq X \) and \( B \in \boldsymbol{\Gamma}(Y) \) there is \( A \in \boldsymbol{\Gamma}(X) \) such that \( B = A \cap Y \). 
A set \(\U \subseteq Y \times X \) is \markdef{$Y$-universal for \(\boldsymbol{\Gamma}(X)\)} if \( \U \in \boldsymbol{\Gamma}(Y \times X) \) and \(\ \boldsymbol{\Gamma}(X)=\left\{ \U_y \mid y \in Y\right \} \), where \(\U_y= \left\{ x \in X \mid (y,x) \in \U \right\} \). It is clear that if \(\U\) is $Y$-universal for \(\boldsymbol{\Gamma}(X)\), then its complement is $Y$-universal for \(\check{\boldsymbol{\Gamma}}(X)\).

\subsection{The \(\I\)-Borel hierarchy}
We recall that a \( \k^+ \)-algebra on set \(X \) is a family of subsets of $X $ closed under the operations of complementation and unions of length at most \( \k\). When $X$ is a topological space, the $\k^+$-algebra generated by the topology of $X$, denoted by \(\Bor(X)\), is the smallest $\k^+$-algebra on $X$ containing all its open sets.

Let \( \nu  \in \{ 2,\k\}\), and consider the topological space \((\pre{\k}{\nu} , \tau_\I) \). We denote the $\k^+$-algebra generated by \(\tau_\I\) by \( \Borid(\pre{\k}{\nu} )\), i.e., \(\Borid(\pre{\k}{\nu} )= \Bor(\pre{\k}{\nu} , \tau_\I)\), and we call its elements \markdef{\( \I \)-Borel} sets. \label{def:Iborel}
When \(\I=b\), \( \Borid(\pre{\k}{\nu} )=\Bor(\pre{\k}{\nu}, \tau_b)\) are the usual generalized Borel sets extensively studied in \cite{MV, FHK14, ACMRP}.

As already observed in \cite{ACMRP}, for every topological space \(X\), the \(\k^+\)-algebra generated by any topology on \(X\) can be stratified in a hierarchy (called the \(\k^+\)-Borel hierarchy) formed by the classes \(\Sii{0}{\a}(X) \), \(\Pii{0}{\a}(X) \), and \( \Dee{0}{\a}(X)  \), where \(\a\) ranges over nonzero ordinals. The pointclasses \(\Bor\), \(\Sii{0}{\a}, \Pii{0}{\a}\) and \(\Dee{0}{\a}\) are boldface, and \(\Bor, \Sii{0}{\a}, \Pii{0}{\a}\) are also hereditary.
In our case study \(X=(\pre{\k}{\nu},\tau_\I)\), the hierarchy of $\I$-Borel sets is defined as follows, and we call it the \markdef{$\I$-Borel hierarchy}.  
 
\begin{definition}\label{def_hierarchy}
The following classes are defined by recursion on the ordinal \(\a\geq 1\):
\begin{align*}
\Siid{0}{1}(\pre{\k}{\nu}) &= \left\lbrace  U \subseteq \pre{\k}{\nu} \mid U \text{ is $\I$-open } \right\rbrace  \\
\Piid{0}{1}(\pre{\k}{\nu})&= \left\lbrace  C \subseteq \pre{\k}{\nu} \mid C \text{ is $\I$-closed } \right\rbrace;  
\end{align*}
\begin{align*}
\Siid{0}{\a}(\pre{\k}{\nu})  &=  \left\lbrace \bigcup_{\gamma < \k} A_\gamma \mid  A_\gamma \in \bigcup_{1 \leq \beta < \alpha}  \Piid{0}{\b}(\pre{\k}{\nu}) \right\rbrace  \\
\Piid{0}{\a}(\pre{\k}{\nu}) &= \left\{ X \setminus A  \mid A \in  \Siid{0}{\a}(\pre{\k}{\nu})\right\}.
\end{align*}
We also set \(  \Deid {0}{\a}(\pre{\k}{\nu}) = \Siid{0}{\a}(\pre{\k}{\nu})  \cap \Piid{0}{\a} (\pre{\k}{\nu}) \).
\end{definition}
 Note that, to be coherent with \cite{ACMRP} and with the usual notation of generalized Borel hierarchy, one should write the classes in Definition~\ref{def_hierarchy} as \(\Sii{0}{\a}(\pre{\k}{\nu}, \tau_\I) \), \(\Pii{0}{\a}(\pre{\k}{\nu}, \tau_\I) \), and \( \Dee{0}{\a}(\pre{\k}{\nu}, \tau_\I)  \). However, for simplicity of notation, we decide to adopt the present one.
 When the space is clear from the context, we drop it from all the notation above.

It is easy to check that
\begin{align}\label{hierarchy_union}
\Borid(\pre{\k}{\nu}) = \bigcup_{1 \leq \alpha < \k^+} \Siid{0}{\a}(\pre{\k}{\nu}) = \bigcup_{1 \leq \alpha < \k^+} \Piid{0}{\a}(\pre{\k}{\nu})   = \bigcup_{1 \leq \alpha < \k^+} \Deid{0}{\a}(\pre{\k}{\nu}).
\end{align}

Since the ideal topology \(\tau_\I\) is finer than the bounded topology \(\tau_b\), a straightforward induction shows the following relations between the classes in the \(\I\)-Borel and the \(\k^+\)-Borel hierarchy.
\begin{proposition}\label{prop: relatioship1} 
For every \(\a<\k^+\), \(\Sii{0}{\a}(\pre{\k}{\nu}, \tau_b) \subseteq \Siid{0}{\a}(\pre{\k}{\nu})\),  hence \\ \(\Pii{0}{\a}(\pre{\k}{\nu}, \tau_b) \subseteq \Piid{0}{\a}(\pre{\k}{\nu}) \) and \(\Dee{0}{\a}(\pre{\k}{\nu}, \tau_b) \subseteq \Deid{0}{\a}(\pre{\k}{\nu})\). Therefore,\\ \(\Bor(\pre{\k}{\nu}, \tau_b) \subseteq \Borid(\pre{\k}{\nu})\). 
\end{proposition}
\begin{proof}
We show that \(\Sii{0}{\a}(\pre{\k}{\nu}, \tau_b) \subseteq \Siid{0}{\a}(\pre{\k}{\nu})\) by induction on \(\a<\k^+\). Since \(b\subseteq \I\), \(\tau_b \subseteq \tau_\I\). Assume \(\a>1\) and that \(\Sii{0}{\b}(\pre{\k}{\nu}, \tau_b) \subseteq \Siid{0}{\b}(\pre{\k}{\nu})\) for every \(\b<\a\). If \(A \in \Sii{0}{\a}(\pre{\k}{\nu}, \tau_b)\), then \(A=\bigcup_{i<\k}A_i\) with \(A_i \in \bigcup_{\b<\a}\Pii{0}{\b}(\pre{\k}{\nu}, \tau_b)\). Since \( \Pii{0}{\b}(\pre{\k}{\nu}, \tau_b) \subseteq \Piid{0}{\b}(\pre{\k}{\nu})\), by the induction hypothesis, \(A_i \in \bigcup_{\b<\a}\Piid{0}{\b}(\pre{\k}{\nu})\) for every \(i<\k\). Hence \(A \in \Siid{0}{\a}(\pre{\k}{\nu})\).
\end{proof}

Recall that \(\B_\I\) denotes the canonical basis for \(\tau_\I\). The following remark shows that \(\B_\I(\pre{\k}{\nu}) \subseteq \Pii{0}{1}(\pre{\k}{\nu}, \tau_b)\).
\begin{remark}\label{rmk: relatioship2} Let \(f \in \Fn_{\I}\) and let \(D =\dom(f)\). 
Then, \(\clopen{f}=\bigcap_{\a \in D} \N^f_\a\), where \(\N^f_\a=\{x \in \pre{\k}{\nu} \mid x(\a)=f(\a)\}\). 
Since \(\N^f_\a \in \Dee{0}{1}(\pre{\k}{\nu}, \tau_b)\) and \(|D|\leq \k\), \(\clopen{f} \in \Pii{0}{1}(\pre{\k}{\nu}, \tau_b)\). 
Thus for every \(f \in \Fn_{\I}\),  \(\clopen{f} \in \Pii{0}{1}(\pre{\k}{\nu}, \tau_b)\).
\end{remark}

\subsection{When \(\I\)-Borel hierarchy is increasing}
 We say that the $\I$-Borel hierarchy is \markdef{increasing} if for every \( 1 \leq \alpha < \beta \) we have \( \Siid{0}{\a}(\pre{\k}{\nu})  \subseteq \Siid{0}{\b}(\pre{\k}{\nu})\).
In this respect, the only problematic case is when \(\a=1\) and \(\b=2\), as 
already noticed more in general in~\cite[Lemma 2.2]{AMR} for any topological space \((X,\tau)\).
Indeed, \( \Siid{0}{\a}(\pre{\k}{\nu})  \subseteq \Siid{0}{\b}(\pre{\k}{\nu})\) for every \( 2 \leq \alpha < \beta \), and  \( \Siid{0}{1}(\pre{\k}{\nu})  \subseteq \Piid{0}{2}(\pre{\k}{\nu}) \subseteq \Siid{0}{3}(\pre{\k}{\nu}) \). In particular, the $\I$-Borel hierarchy is increasing if and only if \(\Siid{0}{1}(\pre{\k}{\nu}) \subseteq \Siid{0}{2}(\pre{\k}{\nu})\).

It was already observed in \cite{HKSW22} that \(\Siid{0}{1}(\pre{\k}{\nu}) \subseteq \Siid{0}{2}(\pre{\k}{\nu})\) when \(\I\) has a basis of size \(\k\).
\begin{proposition}\label{prop: open are kappa unions clopen}\cite[Proposition 1.5]{HKSW22}
     Assume that \(\I\) has a basis of size \(\kappa\). Then every \(\I\)-open subset of \(\pre{\k}{\nu}\) is the union of \(\k\)-many \(\I\)-clopen sets, therefore \(\Siid{0}{1}(\pre{\k}{\nu})\subseteq \Siid{0}{2}(\pre{\k}{\nu})\).
\end{proposition}
\begin{proof}
Let \(\B=\{D_i \mid i<\k\}\) be a basis for \(\I\) and let \(U\subseteq \pre{\k}{\nu}\) be \(\I\)-open. For every \(x \in U\), there exists \(B_x \in \B\) such that \(\clopen{x \restriction B_x}\subseteq U\). Clearly \(U=\bigcup_{x \in U}\clopen{x \restriction B_x} \). Now, for every \(i<\k\), let 
\[
O_i= \bigcup\{\clopen{x \restriction B_x} \mid x \in U \wedge B_x=D_i \}.
\]
If $y\notin O_i$, then $\clopen{y \restriction D_i}\cap O_i=\emptyset$. Thus \(O_i\) is \(\I\)-clopen and \(U=\bigcup_{i<\k}O_i\).
\end{proof}
\begin{corollary}\label{thrm:increasing}
     Assume that \(\I\) has a basis of size \(\kappa\). Then, the $\I$-Borel hierarchy on $\pre{\k}{\nu}$ is increasing.
\end{corollary}

As soon as we remove the assumption on the size of the basis of \(\I\), the $\I$-Borel hierarchy may not be increasing. For example, if \(\I=\sf{NS_\k}\) and we consider the space \(\pre{\k}{2}\), the set \(ub_\k=\{x \in \pre{\k}{2} \mid x=\chi_A{}\footnote{For any $A \subseteq \k$, \(\chi_A \colon \k \to 2\) is the function satisfying \(\chi_A(\a)=1 \Leftrightarrow \a \in A\).}\text{ for some $A \subseteq \k$ unbounded}\}\) is \(\I\)-open but not \( \Siid{0}{2}(\pre{\k}{2})\) (see \cite[Corollary 3.9 and Theorem 3.10]{HKSW22}). 

Note however that certain inclusions hold regardless of whether the \(\I\)-Borel hierarchy is increasing or not. In particular, for all ordinals \(\a\leq \b\) we have 
 \(\Siid{0}{\alpha}(\pre{\k}{\nu})\cup \Piid{0}{\alpha}(\pre{\k}{\nu})\subseteq \Siid{0}{\beta}(\pre{\k}{\nu})\cup \Piid{0}{\beta}(\pre{\k}{\nu})\) and \(\Deid{0}{\a}(\pre{\k}{\nu}) \subseteq \Deid{0}{\b}(\pre{\k}{\nu})\).
.

\begin{remark} \label{rmrk:tall} The set \(ub_\k \subseteq \pre{\k}{2}\) is \(\I\)-open if and only if \(\I\) is tall by \cite[Corollary 3.9]{HKSW22}. Since \(ub_\k \notin \Siid{0}{2}(\pre{\k}{2})\) for every ideal \(\I\) by \cite[Theorem 3.10]{HKSW22}, Proposition~\ref{prop: open are kappa unions clopen} implies that no ideal \(\I\) can simultaneously be tall and have a basis of cardinality \(\k\). 
\end{remark}
Remark \ref{rmrk:tall} implies that any ideal possessing the approximation property cannot be tall. However, the converse implication need not to hold. The following example provides an ideal that is not tall and nevertheless fails to satisfy the approximation property.

\begin{example}\label{messy_counter_example}
   Let $S,T$ be disjoint stationary subsets of \(\kappa\). We define $$\J=\{X\subseteq \kappa \mid X\in \mathsf{NS_\k} \text{ and } X\cap S \text{ is bounded in }\k\}.$$

    Clearly $\J$ is a $\k$-complete ideal that extends the bounded ideal. Notice that $\J$ is not tall, otherwise there would be an unbounded $X\in\J$ such that $X\subseteq S$, contradicting that $X\cap S$ is bounded.

    We now prove that \(\J\) fails to have the approximation property. Assume, for a contradiction, that \(\J\) does have the approximation property, and let \(\langle B_\alpha \mid \alpha<\kappa\rangle\) be an approximation sequence for \(\J\).
   Notice that, for every  $\alpha<\kappa$, $B_\alpha$ is not stationary, so we can choose a club $C_\alpha$ such that $C_\alpha\cap B_\alpha=\emptyset$.
    We inductively construct the sequences $\langle c_\alpha\mid \alpha<\kappa\rangle$ and $\langle t_\alpha\mid \alpha<\kappa\rangle$ such that \(c_\a<t_\a<c_{\a+1}\) and $t_\alpha\in T\cap C_\alpha$. 
    Set \(c_0=0 \), and since $C_0$ is a club and $T$ is a stationary set, pick $t_0 \in T\cap C_0$ such that $c_0<t_0$. 
    Next, suppose $\alpha<\kappa$ is such that for all $\beta<\alpha$ we have chosen $c_\beta$ and $t_\beta$ as desired. If $\alpha$ is a limit ordinal, we let $c_\alpha=\bigcup_{\beta<\alpha}c_\b$; if $\alpha=\b' +1$, is a successor ordinal,  we choose an arbitrary $c_\alpha>t_{\b'}$. Since $C_\alpha$ is a club and $T$ is a stationary set, we choose $t_\alpha\in T\cap C_\alpha$ such that $c_\alpha<t_\alpha$. \\
    Define \(B=\{t_\alpha\mid \alpha<\kappa\}\) and \(C=\{c_\alpha\mid \alpha<\kappa\}\). By construction, \(C\) is club and \(B\cap C=\emptyset\). Furthermore, since \(B\subseteq T\) and \(T\cap S=\emptyset\), it follows that \(B\cap S=\emptyset\). Hence \(B\in \J\). On the other hand, for each \(\alpha<\kappa\) we have \(t_\alpha\in B\cap C_\alpha\) and \(C_\alpha\cap B_\alpha=\emptyset\). Consequently, \(t_\alpha\notin B_\alpha\), and therefore \(B\not\subseteq B_\alpha\) for every \(\alpha<\kappa\). This contradicts the defining property of an approximation sequence for \(\J\). We conclude that \(\J\) does not have the approximation property.
\end{example}

\subsection{When the \(\I\)-Borel hierarchy does not collapse}

We now turn to the question of whether the \(\I\)-Borel hierarchy collapses, as a fundamental parameter measuring the behavior of the \(\I\)-Borel hierarchy is its length. From \eqref{hierarchy_union}, we know that an upper bound for its length is \(\k^+\), but it remained unknown until now whether the \(\I\)-Borel hierarchy on \(\pre{\k}{\nu}\) could be strictly shorter. 
We say that the \(\I\)-Borel hierarchy on \(\pre{\k}{\nu}\) \markdef{collapses} 
if 
\[
\min\left\{\alpha \in \Ord \mid \Siid{0}{\a}(\pre{\k}{\nu}) =\Borid(\pre{\k}{\nu}) \right\}<\k^+.
\]

In the next results we show that for any ideal \(\I\) with the approximation property, the \(\I\)-Borel hierarchy has length \(\k^+\), that is, \(\Siid{0}{\a}(\pre{\k}{\nu}) \subsetneq \Siid{0}{\b}(\pre{\k}{\nu})\) for all \(2 \leq \a<\b<\k^+\), and therefore it does not collapse. This answers \cite[Question 1]{HKSW22} for ideals with the approximation property.

In the generalized context with \(\I=b\), the non-collapse of the \(\k^+\)-Borel hierarchy on $(\pre{\kappa}{2}, \tau_b)$ was first established in \cite[Proposition 4.19]{AMR} for arbitrary infinite cardinals \(\k\), without assuming \(2^{<\k}=\k\). As shown in Theorem~\ref{thrm:collapse_higher_cantor}, this implies that the \(\k^+\)-Borel hierarchy does not collapse on any space containing a homeomorphic copy of $(\pre{\kappa}{2}, \tau_b)$. This result is well-known, but we include the proof for the reader's convenience.

\begin{theorem}\label{thrm:collapse_higher_cantor}
For any topological space $X$ containing an homeomorphic copy of $(\pre{\kappa}{2}, \tau_b)$, the $\kappa^+$-Borel hierarchy on $X$ does not collapse.
\end{theorem}
\begin{proof}
Let $f:\pre{\kappa}{2}\to X$ be a topological embedding and \(Y=\ran(f)\subseteq X \). 
Suppose, towards a contradiction, that \(   \Sii{0}{\a}(X)= \Bor(X)\) for some \(\a<\k^+\).
Since the pointclasses \(\Bor\) and \( \Sii{0}{\a}\) are boldface and hereditary,
    \begin{align*}
    \Sii{0}{\a}(Y)&=\{A \cap Y \mid A \in \Sii{0}{\a}(X)\},\\
     \Bor(Y)&=\{ A \cap Y \mid A \in \Bor(X)\}.
    \end{align*}
Thus, we get \(   \Sii{0}{\a}(Y)= \Bor(Y)\) and the $\kappa^+$-Borel hierarchy on $(\pre{\kappa}{2}, \tau_b)$ collapses, a contradiction. 
\end{proof}

Before stating the main theorem of the section, we present a characterization of the tallness of \(\I\) that will be instrumental in establishing the non-collapse of the \(\I\)-Borel hierarchy under the assumption that \(\I\) is tall.

We first state a result due to Peter Holy, which he conveyed to the authors via private communication.

\begin{fact}[Holy]\label{fact:holy}
 Let $\I=\sf{NS_\k}$, and assume that $2^\mu<2^\kappa$ whenever $\mu<\kappa$. Then, for any $A\subseteq\pre{\k}{2}$, \(A\) endowed with the (induced) \(\I\)-topology is not homeomorphic to $\pre{\k}{2}$ endowed with the bounded topology.
\end{fact}

The next result is the generalization of Fact~\ref{fact:holy}. 

\begin{lemma}\label{lemma:tall_and_homeo2}
 Assume that $\I$ is tall. Then,  for any $A\subseteq\pre{\k}{2}$, \(A\) endowed with the (induced) \(\I\)-topology is not homeomorphic to $\pre{\k}{2}$ endowed with the bounded topology.
\end{lemma}
\begin{proof}
    Suppose towards a contradiction that $\varphi\colon \pre{\k}{2}\to A $ is  an homeomorphism like in the statement. Fix $x\in\pre{\k}{2}$, and for all $\alpha<\kappa$, define 
    $$x_\alpha(\beta)=\begin{cases} 1+_2 x(\beta) &\mbox{if }\beta=\alpha;\\
    x(\beta) & 
    \text{otherwise,}\end{cases}$$ 
    where $+_2$ is the sum in $\mathbb{Z}_2$.

    Let  $y=\varphi(x)$, and for every $\alpha<\kappa$, let $y_\alpha=\varphi(x_\alpha)$. Since $\varphi$ is a bijection, for all $\alpha<\kappa$ we can define $\theta_\alpha<\kappa$ as the smallest ordinal such that $y_\alpha(\theta_\alpha)\neq y(\theta_\alpha)$. Let $S=\{\theta_\alpha\mid \alpha<\kappa\}$.
    \begin{claim}
        $S$ is unbounded in \(\k\).
    \end{claim}
    \begin{proof}[Proof of the Claim.]
        Fix $\gamma<\kappa$. 
        Since $\varphi$ is continuous and $\varphi(x)=y$, there is $\beta<\kappa$ such that $\varphi(\clopen{x\restriction \beta})\subseteq A \cap \clopen{y\restriction\gamma} $. Thus, for all $\alpha>\beta$, $x_\alpha\in \clopen{x\restriction\beta}$, and $y_\alpha\in \clopen{y\restriction\gamma}$. We conclude that $\theta_\alpha>\gamma$.
    \end{proof}
    Since $S$ is unbounded and $\I$ is tall, there exists an  unbounded $D\in \I$ such that $D\subseteq S$. Using again that $\varphi$ is continuous, there exists $\b<\kappa$ such that $\varphi(\clopen{x\restriction \b})\subseteq A \cap \clopen{y\restriction D} $, and thus, for all $\g>\beta$, $y_\g\in A\cap \clopen{y\restriction D}$. 

    Since $\kappa$ is regular, $D\backslash \beta$ is unbounded. Then, there is some $\alpha>\beta$ such that  $\theta_\alpha\in D\backslash \beta$. From the definition of $\theta_\alpha$, $y_\alpha(\theta_\alpha)\neq y(\theta_\alpha)$, so $y_\alpha\notin A\cap \clopen{y\restriction D}$, contradicting the conclusion in previous paragraph.
\end{proof}

The converse of Lemma~\ref{lemma:tall_and_homeo2} holds as well. 

\begin{lemma}\label{lemma:tall_and_homeo1}
    Assume that \(\I\) is not tall. Then, there exists \(C \subseteq \pre{\k}{2}\) such that \((C, \tau_\I\restriction C)\) is homeomorphic to \((\pre{\k}{2}, \tau_b)\).
\end{lemma}
\begin{proof}
Since $\I$ is not tall, there is an unbounded set $S\subseteq \kappa$, such that for all $X\in \I$, $S\cap X$ is bounded.

    Let $$C=\{x\in \pre{\kappa}{2}\mid x(\alpha)=0\text{ for all }\alpha\notin S\}.$$
    Since $\kappa$ is regular, $S$ has order type $\kappa$. Let $e:\kappa\rightarrow S$ be the increasing enumeration of $S$ and define $\Phi: \pre{\k}{2}\rightarrow C$ define by 
    $$\Phi(x)(\alpha)=\begin{cases} x(e^{-1}(\alpha)) &\mbox{if }\alpha\in S;\\
    0 & 
    \text{otherwise.}\end{cases}$$  
    Clearly, $\Phi$ is a bijection. We want to show that $\Phi \colon (\pre{\k}{2}, \tau_b) \rightarrow (C, \tau_\I \restriction C)$ is an homemorphism. 
    
    First, we prove that \(\Phi\) is continuous. Let $y\in C$, $D\in \I$, and $x\in\pre{\kappa}{2}$ such that $\Phi(x)\in C\cap \clopen{y\restriction D}$. Since $D\in\I$, $D\cap S$ is bounded and we get that $e^{-1}(D\cap S)$ is also bounded by some $\theta<\kappa$.  Therefore, $\Phi(z)\in C\cap \clopen{y\restriction D}$ for all $z \in \clopen{x\restriction \theta}$.

    Now we show that \(\Phi\) is open. Let \(x \in \pre{\k}{2}\) and \(\theta<\k\). Notice that $e[\theta]=\{e(\b) \mid \b<\theta\} $ is an increasing sequence of length $\theta$, and since $\kappa$ is regular, \(e[\theta]\) is bounded in $\kappa$ and therefore $e[\theta]\in \I$. 
    We have that \(\Phi(\clopen{x \restriction \theta})=\clopen{t}\cap C\), where \(t \in \pre{e[\theta]}{2}\) is defined by \(t(\alpha)= x(e^{-1}(\alpha))\) for every \(\alpha\in e[\theta] \subseteq S\). 
    Since \(\Phi(x) \in \clopen{t} \cap C\), we are done.
\end{proof}

Combining Lemmas~\ref{lemma:tall_and_homeo2} and~\ref{lemma:tall_and_homeo1} yields Theorem~\ref{tall_and_homeo}, which provides a characterization of the tallness of \(\I\).

\begin{theorem}\label{tall_and_homeo}
An ideal \(\I\) is not tall if and only if there exists a subset \(C \subseteq \pre{\k}{2}\) such that \((C, \tau_\I\restriction C)\) is homeomorphic to \((\pre{\k}{2}, \tau_b)\).
\end{theorem}

We are now ready to state the main result concerning the length of the \(\I\)-Borel hierarchy in the case where \(\I\) is not tall. 
\begin{theorem}\label{non-collapse} 
 If $\I$ is not tall, then the \(\I\)-Borel hierarchy on \(\pre{\k}{2}\) and on \(\pre{\k}{\k}\) does not collapse, therefore it has length \(\k^+\) in both spaces. 
\end{theorem}
\begin{proof}
First, the fact that the \(\I\)-Borel hierarchy on \(\pre{\k}{2}\) has length \(\k^+\) follows from Theorem~\ref{thrm:collapse_higher_cantor} and Theorem~\ref{tall_and_homeo}, because we assume that \(\I\) is not tall.

Next, consider the \(\I\)-Borel hierarchy on \(\pre{\k}{\k}\). 
Suppose, towards a contradiction, that there exists  \(\a<\k^+\) such that  \(\Siid{0}{\a}(\pre{\k}{\k}) =\Borid(\pre{\k}{\k})\). Given \(\nu \in \{2, \k\}\), we recall that in our notation \(\Siid{0}{\a}(\pre{\k}{\nu})=\Sii{0}{\a}(\pre{\k}{\nu} , \tau_\I)\) and \(\Borid(\pre{\k}{\k})=\Bor(\pre{\k}{\k} , \tau_\I) \). 
Since \((\pre{\k}{2}, \tau_\I)\) is a subspace of \((\pre{\k}{\k}, \tau_\I)\), and the classes \(\Sii{0}{\a}\) and \(\Bor\) are hereditary boldface:
 \begin{align*}
 \Sii{0}{\a}(\pre{\k}{2} , \tau_\I) &= \{A \cap \pre{\k}{2} \mid A \in \Sii{0}{\a}(\pre{\k}{\k} , \tau_\I)\},\\
 \Bor(\pre{\k}{2} , \tau_\I)&=\{A \cap \pre{\k}{2} \mid A \in \Bor(\pre{\k}{\k} , \tau_\I)\}.
  \end{align*} 
  Therefore, \(\Siid{0}{\a}(\pre{\k}{2})=\Sii{0}{\a}(\pre{\k}{2} , \tau_\I) =\Bor(\pre{\k}{2} , \tau_\I)=\Borid(\pre{\k}{2})\). A contradiction, since we proved that the \(\I\)-Borel hierarchy on \(\pre{\k}{2}\) does not collapse.
\end{proof}

It remains an open problem whether the \(\I\)-Borel hierarchy collapses for a tall ideal \(\I\). In light of Theorem~\ref{tall_and_homeo}, the method used to establish the analogous statement for not tall ideals does not apply in this context, and consequently additional techniques are required.
\begin{question}
  What is the length of the \(\I\)-Borel hierarchy when \(\I\) is tall?
\end{question}

\subsection{The class of \(\I\)-Borel sets is never trivial}\label{subsec:nontrivial}
Another relevant point is whether the notion of \(\I\)-Borelness is nontrivial, that is, whether \(\Borid(\pre{\k}{\nu}) \neq \pow(\pre{\k}{\nu})\). In some cases, cardinality considerations are useful, for example if \(\I=b\) and \(\k^{<\k}=\k\) then \(|\Bor(\pre{\k}{\nu})|=2^{(\nu^{<\k})} <|\pow(\pre{\k}{\nu})|\). 
However, this approach does not work when we drop the assumption \(2^{<\k}=\k\), or when dealing with ideal topologies. Indeed, if \(\I\) contains an unbounded subset of \(\k\), by Lemma~\ref{lem:basics}\ref{weight}
\[ 
2^{2^\k}=|\tau_\I| \leq |\Borid(\pre{\k}{\nu})| \leq |\pow(\pre{\k}{\nu})|=2^{2^\k}
\] 
therefore \(|\Borid(\pre{\k}{\nu})|=|\tau_\I|=|\pow(\pre{\k}{\nu})|\).

It was already observed in \cite{HKSW22} that, if \(2^{<\k}=\k\), then for any ideal \(\I\) there exists a subset of \(\pre{\k}{2}\) which is not \(\I\)-Borel (\cite[Observation 3.24]{HKSW22}). The same holds when $\I$ is not stationarily tall and \(2^{<\k}>\k\) (\cite[Corollary 3.23]{HKSW22}). 
In the remainder of this subsection, we establish that the concept of \(\I\)-Borelness remains nontrivial in the other missing cases as well. Observe that Theorem~\ref{thrm:withoutBP} resolves \cite[Question 2]{HKSW22}, since it encompasses the previously unresolved case in which \(2^{<\k}>\k\) and \(\I\) is stationarily tall. Furthermore, the theorem subsumes the other cases considered earlier, yielding a unified argument showing that the family of \(\I\)-Borel sets is nontrivial for every \(\k\)-complete ideal \(\I\) that extends the bounded ideal.

The argument used to answer this question, is an $\I$-topology analogue of Petti's Lemma (see \cite[Theorem 9.9]{Kech}).

\medskip

    Let $\mathbb G$ be the group $(2^\kappa,+)$ defined by  $(x+y)(\alpha)=x(\alpha)+_2 y(\alpha)$, where $+_2$ is the sum in $\mathbb{Z}_2$ (the group of two elements). Notice that $\bar 0_\kappa$ is the identity and all the elements have order 2.
    
For all $x\in \pre{\k}{2}$, we denote $Supp(x)=\{\alpha<\kappa\mid x(\alpha)=1\}$, and for $\mathcal J$ ideal on $\kappa$, we let $H_{\mathcal J}=\{x \in \pre{\k}{2}\mid Supp(x)\in \mathcal J\}$. It is straightforward to verify that $H_{\mathcal{J}}$ is a subgroup of $\mathbb G$.

For every $y\in \pre{\k}{2}$, we define \begin{align}
    T_y \colon \pre{\k}{2} \to \pre{\k}{2} \;\;\;\;\; x \mapsto x+y
\end{align} 
and we denote by $T_1$ the case $T_y$ when $y=\bar 1_\kappa$. It is easy to verify that for all $y\in \pre{\k}{2}$, $T_y$ is an homeomorphism, and $T_y=T_y^{-1}$.

\begin{lemma}\label{division_H}
    If $\mathcal{J}$ is a proper maximal ideal, then $\mathbb G=H_\mathcal{J}\cup T_1(H_\mathcal{J})$, and $H_\mathcal{J}\cap T_1(H_\mathcal{J})=\emptyset$.
\end{lemma}
\begin{proof}
    First, let $x \in \mathbb G \setminus H_\mathcal{J}$, so $Supp(x)\notin \mathcal J$. Since $Supp(x + \bar 1_\kappa)=\kappa\setminus Supp(x)$ and $\mathcal{J}$ is maximal, $Supp( x + \bar1_\kappa )\in \mathcal{J}$ and $x + \bar1_\kappa \in H_\mathcal{J}$. Recalling that $x = (x + \bar 1_\kappa) +\bar1_\kappa \in T_1(H_\mathcal{J})$, we are done. 
    For the second part of the statement, notice that for all $x \in H_\mathcal{J}$, $x + \bar 1_\kappa  \notin  H_\mathcal{J}$, otherwise it would be $\bar1_\kappa\in H_J$ and hence $\kappa\in \mathcal{J}$, against $\mathcal{J}$ being proper.
\end{proof}

\begin{lemma}\label{lemma:withoutBP}
    If $\mathcal{I}\subseteq \mathcal{J}$ is a maximal proper ideal, then $H_\mathcal{J}$ is not $\mathcal{I}$-meager and does not have the $\mathcal{I}$-Baire Property.
\end{lemma}
\begin{proof}
First, assume, towards a contradiction, that $H_\mathcal{J}$ is $\mathcal{I}$-meager. Then, since $T_{1}$ is an homeomorphism, $T_1(H_\mathcal{J})$ is $\mathcal{I}$-meager too. Using Lemma~\ref{division_H}, we get that $H_\mathcal{J}\cup T_1(H_\mathcal{J})=\mathbb G$, hence  $\pre{\k}{2}$ is an $\mathcal{I}$-meager set, in contradiction with Proposition \ref{1.3}. Therefore, $H_\mathcal{J}$ is not $\mathcal{I}$-meager in \(\pre{\k}{2}\).

  Now, suppose towards contradiction that $H_\mathcal{J}$ has the $\mathcal{I}$-Baire property. Then, there is some $\mathcal{I}$-open subsets $A $ such that $H_\mathcal{J}\triangle A$ is \(\I\)-meager. Since $H_\mathcal{J}$ is not \(\I\)-meager, we conclude that $A\neq \emptyset$ and $H_\mathcal{J}$ is $\mathcal{I}$-comeager in $A$. 

Moreover, $H_\mathcal{J}$ is $\mathcal{I}$-dense in \(\pre{\k}{2}\), and thus $H_\mathcal{J}\cap A\neq\emptyset$. To see this, let $f\in \Fn_\I(\pre{\k}{2})$. Since $\mathcal{I}\subseteq \mathcal{J}$, $\dom (f)\in \mathcal{J}$. Then, the function \(f^0\) defined as 
 \[
f^0(\g)= \begin{cases}
  f(\g) & \text{if } \g\in\dom( f)\\
  0 & \text{othewise}
  \end{cases}
\]
is such that $Supp(f^0)\subseteq \dom(f)\in \mathcal{J}$, so $f^0 \in H_\mathcal{J}\cap \clopen{f}\neq \emptyset$.

Fix $h\in H_\mathcal{J}\cap A$. We have that $H_\mathcal{J}= T_h(H_\mathcal{J})$: the direct inclusion follows from the fact that $H_\mathcal{J}$ is a subgroup of \(\mathbb G\) and $h\in H_\mathcal{J}$, while the reverse one holds because for all $x\in H_\mathcal{J}$, $T_h(h+x)=x$.
Thus, $T_h(H_\mathcal{J}\cap A)=H_\mathcal{J} \cap T_h(A)$.
Moreover, since \(T_h\) is an homeomorphism, $H_\mathcal{J}$ is \(\I\)-comeager in $T_h(A)$ and $\bar 0_\kappa\in T_h(A)$. 
Therefore, there exists $D \in \mathcal{I}$ such that $\clopen{\bar 0_D}\subseteq T_h(A)$. We denote \(O_D=\clopen{\bar 0_D}\).
Clearly, for all $x,y \in O_D$, $x+y\in O_D$, thus $O_D$ is a subgroup of $ \mathbb G$.

    \begin{claim}\label{O_open}
        $O_D \subseteq H_\mathcal{J}$.
    \end{claim}
    \begin{proof}[Proof of the Claim.]
        For any $x \in O_D$, $x+O_D=O_D$. Since $H_\mathcal{J}$ is $\mathcal{I}$-comeager in $T_h(A)$ and $O_D\subseteq T_h(A)$, $H_\mathcal{J}$ is $\mathcal{I}$-comeager in $O_D$, and hence $x+H_\mathcal{J}$ is $\mathcal{I}$-comeager in $O_D$ too.
        We have that $O_D\cap H_\mathcal{J}\cap (x+H_\mathcal{J})\neq \emptyset$: if not, one could write $O_D=(O_D\backslash H_\mathcal{J})\cup (O_D\backslash x+H_\mathcal{J})$ and since $H_\mathcal{J}$ and $x+H_\mathcal{J}$ are $\mathcal{I}$-comeager in $O_D$, $O_D$ would be $\mathcal{I}$-meager, against Proposition \ref{1.3}.
        Fix $z\in O_D\cap H_\mathcal{J}\cap (x+H_\mathcal{J})$. Then, there is some $y\in H_\mathcal{J}$ such that $z=x+y$, and we get that $x=z+y$. Since both $z,y\in H_\mathcal{J}$, also $x\in H_\mathcal{J}$. We conclude that $O_D\subseteq H_\mathcal{J}$.
    \end{proof}
Finally, we want to show that $H_\mathcal{J}$ is $\mathcal{I}$-clopen. 
Since $H_\mathcal{J}$ is a subgroup, $O_D\subseteq H_\mathcal{J}$, and $\bar 0_\kappa\in O_D$, we have that for all $y\in H_\mathcal{J}$, $T_y(O_D)=y+ O_D\subseteq H_\mathcal{J}$ and $y\in y+ O_D=T_y(O_D)$. Therefore, $H_\mathcal{J}=\bigcup_{y\in H_\mathcal{J}} T_y(O_D)$ is $\mathcal{I}$-open.

To see that $H_\mathcal{J}$ is $\mathcal{I}$-closed, we will show that for all $y\notin H_\mathcal{J}$, $(y+ O_D)\cap H_\mathcal{J}=\emptyset$.
Suppose towards contradiction that there is $y\in H_\mathcal{J}$, $(y+ O_D)\cap H_\mathcal{J}\neq\emptyset$, so there are $z\in H_\mathcal{J}$ and $x\in O_D$  such that $y+x=z$. By Claim \ref{O_open}, $y=z+x \in H_\mathcal{J}$, a contradiction.  Since $\bar0_\kappa\in O_D$ and $O_D$ is $\mathcal{I}$-open, we conclude that $T_y(O_D)$ is an open neighborhood of $y$, thus $H_\mathcal{J}$ is $\mathcal{I}$-closed.

To conclude, we proved that $H_\mathcal{J}$ is $\mathcal{I}$-dense and $\mathcal{I}$-clopen, hence $H_\mathcal{J}=\pre{\k}{2}$. However, since $\mathcal{J}$ is proper, $\kappa\notin \mathcal{J}$ and therefore $\bar 1_\kappa\notin H_\mathcal{J}$. A contradiction. 
\end{proof}

\begin{theorem}\label{thrm:withoutBP}
    There exists a subset of \(\pre{\k}{2}\) that fails to have the $\I$-Baire property. Consequently, there exists a subset of \(\pre{\k}{2}\) that is not \(\I\)-Borel.
\end{theorem}
\begin{proof}
    By Zorn's Lemma, there exists a maximal proper ideal $\mathcal{J}$ extending $\mathcal{I}$. By Lemma~\ref{lemma:withoutBP}, the set $H_\mathcal{J}$ does not have the $\mathcal{I}$-Baire property; in particular, it is not \(\I\)-Borel.
\end{proof}

\subsection{Other properties of the classes in the \(\I\)-Borel hierarchy}
In classical descriptive set theory, universal sets play a crucial role in proving the non-collapse of the Borel hierarchy on uncountable Polish spaces. When moving to the uncountable setting with \(\k>\omega\) and considering \((\pre{\k}{\nu}, \tau_b)\), the situation is as follows. Under the assumption \(2^{<\k}=\k\), there exist \(\pre{\k}{2}\)-universal sets for both \(\Sii{0}{\a}(\pre{\k}{\nu}, \tau_b)\) and \(\Pii{0}{\a}(\pre{\k}{\nu}, \tau_b)\) for each \(1\leq\a<\k^+\) (see, e.g., \cite{ACMRP}). 
On the other hand, if the condition \(2^{<\k}=\k\) fails, then neither \(\Sii{0}{\a}(\pre{\k}{2}, \tau_b)\) nor \(\Pii{0}{\a}(\pre{\k}{2}, \tau_b)\) admits a \(\pre{\k}{2}\)-universal set, for any \(1\leq\a<\k^+\) (\cite[Corollary 4.16]{AMR}). 
The result below shows that when \(\I\) contains an unbounded subset of \(\k\), there are no universal sets for \(\Siid{0}{\a}(\pre{\k}{\nu})\) nor \(\Piid{0}{\a}(\pre{\k}{\nu})\), regardless of whether \(2^{<\k}=\k\) holds.

\begin{proposition}
Assume that \(\I\) contains an unbounded subset of \(\k\). Then, for \(1 \leq \a <\k^+\), neither \(\Siid{0}{\a}(\pre{\k}{\nu})\) nor \(\Piid{0}{\a}(\pre{\k}{\nu})\) have a \(\pre{\k}{2}\)-universal set.
\end{proposition}
\begin{proof}
By Point~\ref{weight} in Lemma~\ref{lem:basics}, \(\Siid{0}{\a}(\pre{\k}{\nu})\) has size greater than \(2^{2^\kappa}\), while \(\left\lbrace \U_y \mid y \in \pre{\k}{2}\right\rbrace \) has size at most \( 2^\k\) for all \(\U \subseteq \pre{\k}{2}\times \pre{\k}{\nu}\). The result for \(\Piid{0}{\a}(\pre{\k}{\nu})\) follows by taking complements.
\end{proof}

We continue with the closure properties of the classes in the \(\I\)-Borel hierarchy on \(\pre{\k}{\nu}\). 
For the proof of the next result, we refer the reader to \cite{ACMRP}. In that paper, the authors consider (regular Hausdorff) topological spaces of weight at most \(\k\) and work under the assumption \(2^{<\kappa}=\kappa\). However, their proof remains valid without the weight restriction. 
Since the \(\I\)-Borel hierarchy may not satisfy the inclusion \(\Siid{0}{1}(\pre{\k}{\nu}) \subseteq \Siid{0}{2}(\pre{\k}{\nu})\), the only case requiring verification is \(\a=1\): by Lemma~\ref{lem:basics}\ref{add}, the class \(\Siid{0}{1}(\pre{\k}{\nu})\) is closed under arbitrary unions and under intersections of fewer than \(\k\) sets.

\begin{proposition}\label{prop: hierarchy_closure}\cite[Proposition 4.1]{ACMRP}
Let \(\nu \in \{2,\k\}\), and $2^{<\kappa}=\kappa$. Given any \( 1 \leq \a < \k^+ \),  let \( \widehat{\alpha} = \k \) if \( \a \) is a successor ordinal, and \( \widehat{\a} = \cof(\a) \) if \( \a \) is limit. Then:
\begin{enumerate}[label={\upshape (\arabic*)}, leftmargin=2pc]
\item 
 \(\Siid{0}{\a}(\pre{\k}{\nu}) \) is closed under unions of length \(\kappa\) and intersections of size smaller than \( \widehat{\a} \);
      \item 
 \(\Piid{0}{\a}(\pre{\k}{\nu}) \)is closed under intersections of length \(\kappa\) and unions of size smaller than \( \widehat{\a} \);
      \item 
\(\Deid{0}{\a}( \pre{\k}{\nu}) \)  is closed under complements and both unions and intersections of size smaller than \( \widehat{\a} \), that is, \(\Deid{0}{\a}( \pre{\k}{\nu}) \)  is a \( \widehat{\a} \)-algebra.
      \end{enumerate} 
\end{proposition}
We refer the reader to \cite{ACMRP} for a detailed discussion of the optimality of the closure properties presented in \ref{prop: hierarchy_closure}.

\subsection{Alternative approaches to Borel sets}

We conclude this section by pointing out a connection with \cite{HMV}. 
When one drops the assumption \(\kappa^{<\kappa}=\kappa\) (equivalently, that \(\kappa\) is regular and \(2^{<\kappa}=\kappa\)), there are two main approaches in the literature. The first maintains the condition \(2^{<\kappa}=\kappa\) but allows \(\kappa\) to be singular, as in \cite{AMR, DM, ACMRP, MRP25}; this ensures that the spaces still have weight \(\kappa\). 
The second approach, developed in \cite{HMV}, keeps \(\kappa\) regular while allowing \(\kappa^{<\kappa}>\kappa\). This constitutes a different setup, as the spaces no longer have weight \(\kappa\). In that context, the authors identify new connections with the model theory of uncountable structures by formulating an alternative notion of Borel sets, tailored to their setting and better suited for model-theoretic applications.
They begin by defining what they call a \emph{basic $\k$-open set}, namely a set of the form \( [f]=\{ x \in \pre{\k}{\nu} \mid f \subseteq x \} \), where \(f \colon X \rightarrow \nu\) and \(X \subseteq \k\) with \(|X| < \k\). As usual, \(\nu \in \{2,\k\}\).

\begin{definition}\label{def_borel_HMV}
 The class of \emph{$\k$-Borel} sets of \(\pre{\k}{\nu}\)  is the smallest class containg the basic $\k$-open sets and closed under complements and unions of size at most $\k$.
\end{definition}
In line with Definition~\ref{def_borel_HMV}, one could attempt to develop the theory of Borel sets in ideal topologies using the following definition instead of the one adopted in this article (cf. page~\ref{def:Iborel}). 

\begin{definition}\label{def_borel_basic}
    The class of \markdef{basic $\I$-Borel} of \(\pre{\k}{\nu}\) is the smallest class containg the basic $\I$-open sets and closed under complements and unions of size at most $\k$.
\end{definition}

As the next result shows, for any ideal \(\I\) the class of basic \(\I\)-Borel sets coincides with the class defined in Definition~\ref{def_borel_HMV}.

\begin{fact}
    For any ideal $\I$, the class of basic $\I$-Borel sets  of \(\pre{\k}{\nu}\) coincide with the class of $\k$-Borel sets  of \(\pre{\k}{\nu}\). As a consequence, the class of basic $\I$-Borel is independent of the ideal $\I$.
\end{fact}
\begin{proof}
As \(b \subseteq \I\), any basic $\k$-open set is a basic $\I$-open set, therefore every $\k$-Borel set is a basic $\I$-Borel set. 
For the other inclusion, it is sufficient to show that all basic $\I$-open sets are $\k$-Borel sets. As observed in Remark~\ref{rmk: relatioship2}, for every \(f \in \Fn_\I\), \(\clopen{f}=\bigcap_{\a \in D} \N^f_\a\), where  $D =\dom(f)$ and \(\N^f_\a=\{x \in \pre{\k}{\nu} \mid x(\a)=f(\a)\}\). Since each \(\N^f_\a\) is a basic \(\k\)-open sets and \(|D|\leq \k\), \(\clopen{f}\) is \(\k\)-Borel.
\end{proof}

\section{The \(\I\)-Analytic sets}\label{section_analytic}

In this section we introduce the notion of generalized analytic set for the spaces \(\pre{\k}{2}\) and \(\pre{\k}{\k}\) endowed with the ideal topology \(\tau_\I\). However, we immediately observe that this notion is trivial, because the class under consideration actually coincides with the entire powerset (Corollary~\ref{all_analitic}). Moreover, the same phenomenon occurs for the generalized versions of all the equivalent characterizations of analytic sets from classical descriptive set theory (Corollary~\ref{analequivalence}). As usual, fix \(\nu\in \{2,\k\}\).

Given \(\mu\in \{2,\k\}\), a function \(\Phi \colon \pre{\k}{\nu} \to \pre{\k}{\mu}\) 
is \markdef{\(\I\)-Borel} if \(\Phi^{-1}(U) \in \Borid(\pre{\k}{\nu})\) for every \(\I\)-open set \(U\subseteq \pre{\k}{\mu}\).
We say that a set $B \subseteq \pre{\k}{\mu}$ is an \(\I\)-continuous image of $A\subseteq \pre{\k}{\nu}$ if there is an \(\I\)-continuous function \(\Phi \colon \pre{\k}{\nu} \to \pre{\k}{\mu}\) such that $\Phi(A)=B$. If $\Phi$ is \(\I\)-Borel, then we say that $B$ is an \(\I\)-Borel image of $A$.

\begin{definition}\label{def:anal}
    A set \(A \subseteq \pre{\k}{\nu}\) is called \markdef{\(\I\)-analytic} (\(\Siid{1}{1}(\pre{\k}{\nu})\)) if it is either empty or an \(\I\)-continuous image of an \(\I\)-closed subset of \(\pre{\k}{\k}\).
\end{definition}

When \(\I=b\), \( \Siid{1}{1}(\pre{\k}{\nu} )\) is the usual class of \(\k\)-analytic sets as in \cite{FHK14, LS15, AMR}.

\begin{theorem}\label{anal_closure}
If \(\I\) contains an unbounded subset of \(\k\), then the class \(\Siid{1}{1}(\pre{\k}{\nu})\) is closed under unions of size at most \( 2^\kappa\).
\end{theorem}
\begin{proof} 
Let \(\seq{A_\a}{\a<2^\k}\) be a sequence of \(\Siid{1}{1}\)-subsets of \(\pre{\k}{\nu}\). For every \(\a<2^\k\), let \(C_\a\subseteq \pre{\k}{\k}\) be \(\I\)-closed and \(\Phi_\a \colon \pre{\k}{\k} \to \pre{\k}{\nu}\) be \(\I\)-continuous such that \(\Phi_\a(C_\a)=A_\a\). 
Since $\I$ contains an unbounded set, there are \(U_0,U_1 \subseteq \k\) unbounded such that \(U_0,U_1 \in \I\) and \(U_0\cap U_1 =\emptyset\). 
Fix the bijections
\begin{align*}
    &\theta \colon \pre{U_1}{\k} \to 2^\kappa \\
    &\d \colon U_0 \to \k
\end{align*}
and define 
\(\Delta \colon \pre{\k}{\k} \to \pre{\k}{\k}\) such that 
\(\Delta(x)(\a)=x(\d^{-1}(\a))\) for all \(\a<\k\).
Given \(x \in \pre{\k}{\k}\), notice that for every \(y \in \clopen{x \restriction U_0}\), \(\Delta(y)=\Delta(x)\) and, for every \(z \in \clopen{x \restriction U_1}\), \(\theta(z\restriction U_1)=\theta(x\restriction U_1)\).
Finally, let
\[
\Phi \colon \pre{\k}{\k} \to \pre{\k}{\nu} \colon x \mapsto \Phi_{\theta(x\restriction U_1)}(\Delta(x)),
\]
and define 
\[
\; C=\{ x \in \pre{\k}{\k} \mid \Delta(x)\in C_{\theta(x\restriction U_1)}\}.
\]
We claim that \(\Phi\) and \(C\) witness \(\bigcup_{\a<2^\k}A_\a \in \Siid{1}{1}(\pre{\k}{\nu})\). 

Notice that for every \(x \in\pre{\k}{\k} \) and for every \(z \in \clopen{x \restriction U_0 \cup U_1}\),
\(\Phi(z)=\Phi_{\theta(x\restriction U_1)}(\Delta(x))=\Phi(x)\), therefore \(\Phi(\clopen{x \restriction U_0 \cup U_1})=\{\Phi(x)\}\). 

 The function \(\Phi\) is \(\I\)-continuous since for every \(x \in \pre{\k}{\k}\), and for every \(D \in \I\), the set \(\clopen{x \restriction U_0 \cup U_1}\) is \(\I\)-open (because \(U_0\cup U_1 \in \I\)) and it satisfies \(\Phi(\clopen{x \restriction U_0 \cup U_1})\subseteq \clopen{\Phi(x) \restriction D}\).
To see that \(C\) is \(\I\)-closed we show that for every \(x \notin C\), \(\clopen{x \restriction U_0 \cup U_1} \subseteq \pre{\k}{\k}\setminus C\). This is immediate, as for every \(y \in \clopen{x \restriction U_0 \cup U_1} \), \(\Delta(y)=\Delta(x)\) and \(\theta(y\restriction U_1)=\theta(x\restriction U_1)\), therefore \(\Delta(y) \notin C_{\theta(y\restriction U_1)}\). The following claim concludes the proof.
\begin{claim}
\(\Phi(C)=\bigcup_{\a<2^\k}A_\a\).
\end{claim}
\begin{proof}[Proof of the Claim.]
    If \(y \in \bigcup_{\a<2^\k}A_\a\), then there is some \(\a<2^\k\) and some \(x\in C_\a \) such that \(y=\Phi_\a(x)\). 
Let \(z \in \pre{\k}{\k}\) such that \(\theta(z \restriction U_1)=\a\) and \(\Delta(z)=x\).  
Then, 
\[
\Phi(z)=\Phi_{\theta(z\restriction U_1)}(\Delta(z))=
\Phi_{\alpha}(x)=y.
\]
Finally, since $x\in C_\alpha$, $\Delta(z)\in C_{\theta(z\restriction U_1)}$. Thus $z\in C$ and \(y \in \Phi(C)\).

 For the other direction, let \(y \in C\). From the definition of $C$, $\Delta(y)\in C_{\theta(y\restriction U_1)}$. Let \(\a=\theta(y \restriction U_1)\) and \(x=\Delta(y)\). Then, \(x \in C_\a\) and $\Phi_\a(x)\in A_\a$. Therefore,  
\(\Phi(y)=\Phi_{\theta(y\restriction U_1)}(\Delta(y))=\Phi_\a(x)\in A_\a\).
\end{proof}\end{proof} 

\begin{corollary}\label{all_analitic}
    If \(\I\) contains an unbounded subset of \(\k\), then \(\Siid{1}{1}(\pre{\k}{\nu})=\pow(\pre{\k}{\nu})\).
\end{corollary}


The generalized Cantor and Baire spaces (equipped with the bounded topology, hence \(\I=b\)) have \(2^{(2^{<\k})}\)-many \(\k\)-analytic subsets, so the hypothesis \(2^{<\k}=\k\) guarantees that the class \(\Siid{1}{1}\) is strictly smaller than the full powerset. In contrast, when \(\k^{<\k}>\k\), \cite[Lemma 2.1]{HMV} shows that \(\Siid{1}{1}(\pre{\k}{2})=\pow(\pre{\k}{2})\) occurs even if $\I$ is the bounded ideal.

Since in classical descriptive set theory the notion of analytic set can be reformulated giving several equivalent definitions \cite[Section 14.A]{Kech}, the reader may wonder if the generalized counterpart of any of the other definitions may be nontrivial. This is not the case, as the next result shows.

 Given a subset \(A \subseteq X \times Y\) of a cartesian product, we denote its projection (on the first coordinate) by \(\proj(A)=\{x \in X \mid \exists y \in Y (x,y) \in A \}\). 

\begin{corollary}\label{analequivalence}
	If \(\I\) contains an unbounded subset of \(\k\), the following are equivalent:
	\begin{enumerate}[label={\upshape (\arabic*)}, leftmargin=2pc]
     \item \label{analequivalence6} \(A \in \pow(\pre{\k}{\nu})\)
		\item\label{analequivalence1} A is \(\I\)-analytic;
		\item\label{analequivalence3} \(A = \proj(C)\) for some \(\I\)-closed \(C \subseteq \pre{\kappa}{\nu} \times \pre{\k}{\kappa}\);
        \item\label{analequivalence2} A is an \(\I\)-continuous image of an \(\I\)-Borel subset of \( \pre{\k}{\kappa}\);
		\item\label{analequivalence4} \(A = \proj(B)\) for some \(\I\)-Borel \(B \subseteq \pre{\kappa}{\nu} \times \pre{\k}{\kappa}\);
	\item\label{analequivalence5} A is an \(\I\)-Borel image of an \(\I\)-Borel subset of  \( \pre{\k}{\kappa}\).
\end{enumerate}
\end{corollary}
\begin{proof}
Trivially, \ref{analequivalence1}-\ref{analequivalence5} implies \ref{analequivalence6}.
Clearly \ref{analequivalence1} implies \ref{analequivalence2}. \ref{analequivalence3} implies \ref{analequivalence4}, since every \(\I\)-closed set is \(\I\)-Borel. \ref{analequivalence2} implies \ref{analequivalence5}, since every \(\I\)-continuous function is \(\I\)-Borel. By Corollary~\ref{all_analitic}, \ref{analequivalence6} implies \ref{analequivalence1}. 
    It remains to show that \ref{analequivalence1} implies \ref{analequivalence3}. Let \(A\subseteq \pre{\k}{\nu} \) be \(\I\)-analytic and let \(\Phi \colon C \to \pre{\k}{\nu} \) be an \(\I\)-continuous surjection onto \(A\) for some \(\I\)-closed \(C \subseteq \pre{\k}{\k}\). Since
    \((\pre{\k}{\nu}, \tau_\I)\) is a Hausdorff space by Lemma~\ref{lem:basics}\ref{T2}, this implies that \(\operatorname{Graph}(\Phi)=\{(x,y) \in C \times \pre{\k}{\nu} \mid \Phi(x)=y\}\) is  \(\I\)-closed in \(C \times \pre{\k}{\nu}\), hence also \(\I\)-closed in \(\pre{\k}{\k} \times \pre{\k}{\nu} \). Clearly, \(A=\proj(\{(y,x) \in \pre{\k}{\nu} \times \pre{\k}{\k} \mid (x,y) \in \operatorname{Graph}(\Phi)\})\).
    \end{proof}

Only one possible definition of \(\I\)-analytic set is missing from Corollary~\ref{analequivalence}.
In classical descriptive set theory (hence \(\k=\o\)), it is well-known (see \cite[Section 14-15]{Kech}) that a non-empty set \(A \subseteq \pre{\o}{2}\) is:
\begin{itemize}
    \item a continuous image of \(\pre{\o}{\o}\) if and only if \(A\) is analytic;
    \item an injective continuous image of some closed \(C \subseteq \pre{\o}{\o}\) if and only if \(A\) is Borel. 
\end{itemize}

Moreover, in the generalized setting \(\k^{<\k}=\k>\o\), the class of  continuous injective images of closed subsets of \((\pre{\k}{\k}, \tau_b)\) consistently coincide with the class of \(\k\)-analytic sets \cite{LS15}.
It is therefore natural to ask whether the generalization of at least one of these two definitions may give a nontrivial notion of \(\I\)-analytic set. The answer is negative again.  

\begin{proposition} \label{prop_wrong_anal}
Assume that \(\I\) contains an unbounded subset of \(\k\). Then:
	\begin{enumerate}[label={\upshape (\arabic*)}, leftmargin=2pc]
		\item\label{prop_wrong_anal1} Every non-empty \(A \subseteq \pre{\k}{\nu}\) is an \(\I\)-continuous image of \(\ \pre{\k}{\k} \).
		\item\label{prop_wrong_anal2} Every \(A \subseteq \pre{\k}{\nu}\) is an injective \(\I\)-continuous image of some \(\I\)-closed \(C \subseteq \pre{\k}{\k}\).
	\end{enumerate}
\end{proposition}

\begin{proof}
Let \(U\in \I\) be unbounded and fix a bijection \(\theta:\pre{U}{\k} \rightarrow \pre{\k}{\nu}\). Then  
\[
\Phi:\pre{\k}{\k} \rightarrow \pre{\k}{\nu}: x \mapsto \theta(x\restriction {U})
\] 
is such that for every non-empty \(A \subseteq \pre{\k}{\nu}\), \(\Phi^{-1}(A)\) is \(\I\)-clopen, because 
\[
\Phi^{-1}(A)=\bigcup \left\lbrace \Nbhd_f \mid f \in \pre{U}{\k}, \theta(f) \in A \right\rbrace 
\] 
and
\[ 
\pre{\k}{\k} \setminus \Phi^{-1}(A)=\bigcup \left\lbrace \Nbhd_f \mid f \in \pre{U}{\k}, \theta(f) \notin A \right\rbrace 
\] 
are both \(\I\)-open. In particular, \(\Phi:\pre{\k}{\k} \rightarrow \pre{\k}{\nu}\) is \(\I\)-continuous. 

To show \ref{prop_wrong_anal1}, we define the surjective function \(\Psi: \pre{\k}{\k} \to A \) by setting:
 \begin{equation*}
\Psi(x)=
\begin{cases}
\Phi(x) & \text{if \(x \in \Phi^{-1}(A)\)} \\
x_0 & \text{if \(x \notin \Phi^{-1}(A)\)}
\end{cases} 
\end{equation*} for some \(x_0 \in A\). Clearly, \(\Psi\) is continuous. 

For \ref{prop_wrong_anal2}, consider the set 
\[
C=\left\lbrace x \in \pre{\k}{\k} \mid x(\a) = 0 \text{ for all } \a \in \k\setminus U, \theta(x \restriction U) \in A \right\rbrace.
\] 
Since \(|\k\setminus U|\leq \k\) and \(\left\lbrace x \in \pre{\k}{\k} \mid x(\a) \neq 0  \right\rbrace \in \Dee{0}{1} \subseteq \Deid{0}{1}\), \[
\pre{\k}{\k}\setminus C=\bigcup_{\a \in \k\setminus U} \left\lbrace x \in \pre{\k}{\k} \mid x(\a) \neq 0  \right\rbrace \cup \bigcup \left\lbrace \Nbhd_f \mid \theta(f)\notin A \right\rbrace
\]
is \(\I\)-open. Thus \(C\) is \(\I\)-closed.
Then \(\Phi\restriction C\) is injective by construction (since \(\theta\) is bijective) and it is onto \(A\).
\end{proof}

\section{The Approximation Lemma}\label{section_approx}

The Approximation Lemma has been one of the main techniques for constructing continuous reductions between equivalence relations in the generalized Baire space (i.e. endowed with the bounded topology under the assumption $\kappa^{<\kappa}=\kappa$). 
In \cite{FMR21} this method was also extensively used to study reductions between different equivalence relations under some reflection principles. We refer the reader to \cite{FMR21} for all undefined notions in this section. 

One of the main applications of the Approximation Lemma in the bounded topology is the characterization of some reductions, for instance the ones in Example~\ref{example_approx}.
In this framework, if $S,X\subseteq \k$ are disjoint stationary sets, then the following are equivalent:
\begin{enumerate}[label={\upshape (\arabic*)}, leftmargin=2pc]
    \item $=_S^\k$ is Lipschitz reducible to $=_X^\k$; 
    \item $=_S^\k$ is $X$-recursive reducible to $=_X^\k$;
    \item $=_S^\k$ has an $X$-approximation with the club filter\footnote{See Example \ref{example_approx} for the definition of $=_S^\k$, and Definition \ref{A-approximationk} for the definition of $X$-approximation.}.
\end{enumerate}
All these considerations motivate the study of the Approximation Lemma for ideal topologies.
In this section we study the Approximation Lemma for the space \((\pre{\k}{\k}, \tau_\I)\).
Since we consider only ideals that have the approximation property  (Definition \ref{approx_property_for_last}),  throughout this section we assume that $\mathcal{I}$ has a basis of size $\kappa$.

Given an approximation sequence $\langle B_\alpha\rangle_{\alpha<\kappa}$ of $\mathcal{I}$, we define the search function $s:\kappa\rightarrow\kappa$ as follows:
$$s(\alpha)=\begin{cases} \beta &\mbox{if }\alpha\in B_{\beta+1}\backslash B_\beta \mbox{ and } \beta<\alpha,\\
    \alpha & 
    \text{otherwise.}\end{cases}$$  
Note that $s$ is well-defined and surjective, since $\langle B_\alpha\rangle_{\alpha<\kappa}$ is an strictly increasing continuous sequence and $B_0=\emptyset$.
Notice that the search function tell us whether an ordinal $\alpha$ appears in the approximation sequence before $B_{\alpha+1}$ and in case it does, it tells us the maximum $B_\beta$ that does not contains $\alpha$. 

\begin{definition}
    Let $\seq{B_\alpha}{\alpha<\kappa}$ be an approximation sequence of $\mathcal{I}$, and $A\subseteq \kappa$ an unbounded set. We say that a function $f:\pre{\k}{\kappa}\rightarrow\pre{\k}{\kappa}$ is \markdef{$A$-recursive over $\langle B_\alpha\rangle_{\alpha<\kappa}$} if there exists a monotone function $$H:\Fn_{\I}\rightarrow \Fn_{\I}$$ such that
    \begin{itemize}
        \item for every $x\in \pre{\k}{\kappa}$ and $\beta<\kappa$, $\dom (H(x\restriction B_\beta))=s^{-1}[\beta+1]$.
        \item for every $\alpha<\kappa$ and for every $x\in \pre{\k}{\kappa}$, 
           $$f(x)(\alpha)=H(x\restriction B_{s(\alpha)})(\alpha).$$ 
        \item for every $\alpha\notin A$ and for every $x\in \pre{\k}{\kappa}$, $f(x)(\alpha)=0$.
    \end{itemize}
           
\end{definition}

Notice that the role of $A$ is to be the set in which $f$ may be not zero. We may define a recursive function over $\langle B_\alpha\rangle_{\alpha<\kappa}$ by making $A=\kappa$. From the definition of $s$, for all $\alpha\in B_{\beta+1}$, we have that $s(\alpha)\leq \beta$. Thus $B_{\beta+1}\subset s^{-1}[\beta+1]$. 

Let $D\in \I$, since $\seq{B_\alpha}{\alpha<\kappa}$ is an approximation sequence of $\mathcal{I}$, there is $\beta<\kappa$ such that $D\subseteq B_\beta$. Since $B_{\beta}\subseteq \dom (H(x\restriction B_\beta))$, for all $x\in \pre{\k}{\kappa}$, $D\subseteq \dom (H(x\restriction B_\alpha))$.
Therefore, if \(H\) witnesses that $f$ is $A$-recursive over $\seq{B_\alpha}{\alpha<\kappa}$, then $H$ is continuous as in Definition \ref{def_functions}. In particular \(H^*=f\), 
thus \(f\) is \(\I\)-continuous by Proposition \ref{prop:cont}. 

\begin{fact}
    Let $\seq{B_\alpha}{\alpha<\kappa}$ be an approximation sequence of $\mathcal{I}$, and $A\subseteq \kappa$ an unbounded set. If $f$ is $A$-recursive over $\seq{B_\alpha}{\alpha<\kappa}$, then $f$ is \(\I\)-continuous.
\end{fact}

Let $E_1$ and $E_2$ be equivalence relations on $\pre{\k}{\kappa}$. Let $A\subseteq \kappa$ be an unbounded set.
We say that \markdef{$E_1$ is $A$-recursive reducible to $E_2$} if there are an approximation sequence $\seq{B_\alpha}{\alpha<\kappa}$ of $\mathcal{I}$, and a function $f\colon
\pre{\k}{\kappa}\rightarrow \pre{\k}{\kappa}$ $A$-recursive over $\seq{B_\alpha}{\alpha<\kappa}$ that satisfies $$(x,y)\in E_1 \Leftrightarrow
(f(x),f(y))\in E_2.$$  We call $f$ an \markdef{$A$-recursive reduction} of $E_1$ to
$E_2$, and when such reduction exists we write $E_1\hookrightarrow_{R(A)} E_2$. Clearly, $E_1\hookrightarrow_{R(A)} E_2$ implies that an \(\I\)-continuous (\(\I\)-Borel) reduction exists.
Thus an $A$-recursive reduction is stronger than an \(\I\)-continuous reduction.

\begin{definition}\label{def:filter_relation}
    Let $\mathcal{F}\subseteq \mathcal{P}(\kappa)$ be a filter. 
    We define the following equivalence relation $\sim_\mathcal{F}^{\kappa}$ on the space $\pre{\k}{\kappa}$: 
    $$x \sim^{\kappa}_\mathcal{F}y\Leftrightarrow \exists W\in\mathcal{F}\ \left( W\subseteq \{\alpha<\kappa\mid x(\alpha)=y(\alpha)\}\right).$$
\end{definition}

Given a filter $\mathcal{F}$, we denote by $\mathcal{F}^+$ the set $\{X\subseteq \kappa\mid \forall Y\in \mathcal{F} (X\cap Y\neq \emptyset)\}$. For every $A\in \mathcal{F}^+$, we define the filter $\mathcal{F}\restriction A$ as the filter generated by $\{X\cap A\mid X\in \mathcal{F}\}$.

\begin{example}\label{example_approx}
Many well-known equivalence relation are induced by a filter like in Definition~\ref{def:filter_relation}. The following are some examples. 
\begin{itemize}
    \item The identity relation $\id$ is the equivalence  relation induced by the filter $\{\k\}$. 
    \item For \(X \subseteq \k\), the equivalence modulo $X$, denoted by $X_\k$, is the equivalence relation induced by the filter $\{\k\}\restriction X=\{Y\subseteq \k\mid X\subseteq Y\}$. 
    \item $E_0^{<\kappa,\kappa}=\{(x,y)\in \pre{\k}{\kappa}\times \pre{\k}{\kappa}\mid \exists\alpha<\kappa \forall\beta>\alpha\ (x(\beta)=y(\beta))\}$ is the equivalence relation induced by the filter generated by $\{\k\setminus \alpha \mid \alpha<\k\}$.
    \item The equivalence modulo nonstationary set $=_{\CUB}^\k$ is the equivalence relation induced by the club filter $\CUB$.
    \item For $S$ be a stationary subset of \(\k\), $=_S^\k$ is the equivalence relation induced by the filter $\CUB\restriction S$.
\end{itemize}
\end{example}

\begin{definition}\label{A-approximationk}
    Let \(\seq{B_\alpha}{\alpha<\kappa}\) be an approximation sequence of $\mathcal{I}$, and let  $\mathcal{H}$ be a filter over $\kappa$ with basis $\mathcal{B}$. Let $A\in \mathcal{H}^+$ be an unbounded set and let $E$ be an equivalence relation on $\pre{\k}{\kappa}$.
    We say that $E$ has an \markdef{$A$-approximation over $\seq{B_\alpha}{\alpha<\kappa}$} if there is  a sequence $\langle E_\alpha\mid \alpha<\kappa\rangle$ with $E_\alpha\subseteq \pre{B_{s(\alpha)}}{\k}\times \pre{B_{s(\alpha)}}{\k}$ such that the following conditions hold:
    \begin{enumerate}[label={\upshape (\arabic*)}, leftmargin=2pc]
        \item For every $ \alpha<\kappa$, $E_\alpha$ is an equivalence relation.
        \item\label{A-approximationk_2} For every $x,y\in \pre{\k}{\kappa}$, if $x E y$ then there exists $D\in \mathcal{B}$ such that for every $\alpha\in D$, $x\restriction B_{s(\alpha)} \ E_\alpha\ y\restriction B_{s(\alpha)}$.
        \item\label{A-approximationk_3} For every $x,y\in \pre{\k}{\kappa}$, if $\neg (x E y)$ then there exists $A'\in \mathcal{H}^+$, $A'\subseteq A$, such that for every $\alpha\in A'$, $\neg(x\restriction B_{s(\alpha)} \ E_\alpha\ y\restriction B_{s(\alpha)})$.
    \end{enumerate}
\end{definition}
An $A$-approximation is said to be \markdef{precise} if for every $\alpha<\kappa$, $E_\alpha$ has at most $\kappa$-many classes. Clearly, the assumption $\k^{<\k}=\k$ implies that every $A$-approximation is precise.
We say that $E$ has an (precise) $A$-approximation if there is an approximation sequence $\seq{B_\alpha}{\alpha<\kappa}$ of $\mathcal{I}$, such that $E$ has an (precise) $A$-approximation over $\seq{B_\alpha}{\alpha<\kappa}$.


\begin{lemma}[Approximation Lemma]\label{approx_lemma}
    Let $\mathcal{H}$ be a filter over $\kappa$ with basis $\mathcal{B}$. Let $A\in \mathcal{H}^+$ be an unbounded set and let $E$ be an equivalence relation on $\pre{\k}{\kappa}$. Then, the following are equivalent:
    \begin{enumerate}[label={\upshape (\arabic*)}, leftmargin=2pc]
        \item\label{approx1} $E$ has a precise $A$-approximation;
        \item\label{approx2} $E\hookrightarrow_{R(A)}\sim^\kappa_\mathcal{F}$, where $\mathcal{F}=\mathcal{H}\restriction A$.
\end{enumerate}
\end{lemma}

\begin{proof}
\ref{approx1}$\Rightarrow$ \ref{approx2}.
Let $\seq{B_\alpha}{\alpha<\kappa}$ be an approximation sequence of $\mathcal{I}$ such that there is $\langle E_\alpha\mid \alpha<\kappa\rangle$ a precise $A$-approximation over $\seq{B_\alpha}{\alpha<\kappa}$. For every $\alpha<\kappa$, let $\langle x^\alpha_i\mid 0<i<\kappa\rangle$ be an enumeration of the equivalence classe $E_\alpha$. 
    Let us define $F:\pre{\k}{\kappa}\rightarrow\pre{\k}{\kappa}$ as follows:
        $$F(x)(\alpha)=\begin{cases} i &\mbox{if }\alpha\in A \mbox{ and } x\restriction B_{s(\alpha)}\in x_i^{\alpha},\\
    0 & \mbox{otherwise. }\end{cases}$$  
    Let us show that $x\ E\ y$ if and only if $F(x)\ \sim^\kappa_\mathcal{F}\ F(y)$.
    \begin{claim}
        $x\ E\ y$ implies $F(x)\ \sim^\kappa_\mathcal{F}\ F(y)$.
    \end{claim}
    \begin{proof}[Proof of the Claim.]
        Suppose $x,y\in \pre{\k}{\kappa}$ are such that $x\ E\ y$. Since $\langle E_\alpha\mid \alpha<\kappa\rangle$ is an $A$-approximation over $\seq{B_\alpha}{\alpha<\kappa}$, by Definition~\ref{A-approximationk}\ref{A-approximationk_2} there exists $D\in\mathcal{B}$ such that for every $\alpha\in D$, $$x\restriction B_{s(\alpha)}\ E_\alpha\ y\restriction B_{s(\alpha)}.$$ So, for all $\alpha\in D\cap A$, $F(x)(\alpha)=F(y)(\alpha)$. Since $D\cap A\in \mathcal{H}\restriction A=\mathcal{F}$, we get $F(x)\ \sim^\kappa_\mathcal{F}\ F(y)$.
    \end{proof}
    \begin{claim}
        $\neg(x\ E\ y)$ implies $\neg(F(x)\ \sim^\kappa_\mathcal{F}\ F(y))$.
    \end{claim}
    \begin{proof}[Proof of the Claim.]
        Suppose $x,y\in \pre{\k}{\kappa}$ are such that $\neg(x\ E\ y)$. Since $\langle E_\alpha\mid \alpha<\kappa\rangle$ is an $A$-approximation over $\seq{B_\alpha}{\alpha<\kappa}$, by Definition~\ref{A-approximationk}\ref{A-approximationk_3}, there is $A'\subseteq A$ in $\mathcal{H}^+$,  such that for all $\alpha\in A'$, $$\neg(x\restriction B_{s(\alpha)}\ E_\alpha\ y\restriction B_{s(\alpha)}).$$ So, for all $\alpha\in A'$, $F(x)(\alpha)\neq F(y)(\alpha)$. Since $A'\in\mathcal{H}^+$, for all $X\in \mathcal{H}$, $X\cap A'\neq\emptyset$. Thus, for every $X\in \mathcal{H}$, $(X\cap A)\cap A'\neq\emptyset$. So $A'\in (\mathcal{H}\restriction A)^+=\mathcal{F}^+$, and
        we conclude that $\neg(F(x)\ \sim^\kappa_\mathcal{F}\ F(y))$.
    \end{proof}
    \begin{claim}
        $F$ is $A$-recursive over $\seq{B_\alpha}{\alpha<\kappa}$.
    \end{claim}
    \begin{proof}[Proof of the Claim.]
    Define $H:\Fn_{\I}\rightarrow \Fn_{\I}$ as follows:
        $$H(x\restriction D)=F(x)\restriction s^{-1}[\alpha+1]$$
    where $\alpha$ is the smallest element such that $D\subseteq B_\alpha$.
    Notice that if $D=B_\alpha$, then $H(x\restriction D)=H(x\restriction B_\alpha)= F(x)\restriction s^{-1}[\alpha+1]$. Also, for all $\beta<\kappa$, $\beta\in s^{-1}[s(\beta)+1]$ and $H(x\restriction B_{s(\beta)})(\beta)=F(x)(\beta)$. Finally $F$ is 0 outside $A$. 

   Thus we are only missing to show that $H$ is well defined, it is sufficient to check the case $\dom(x)=B_\alpha$. Let $x,y\in \kappa^{\kappa}$ be such that $x\restriction B_\alpha=y\restriction B_\alpha$. 
    Clearly, for every $\beta\in s^{-1}[\alpha+1]\backslash A$, $F(x)(\beta)=0=F(y)(\beta)$, therefore $F(x)\restriction s^{-1}[\alpha+1]\ (\beta)=0=F(y)\restriction s^{-1}[\alpha+1]\ (\beta)$ for all $\beta\in s^{-1}[\alpha+1]\backslash A$.
    On the other hand, let $\beta\in s^{-1}[\alpha+1]\cap A$. From the definition of $F$, $F(x)(\beta)=i$ and $F(y)(\beta)=j$, where $x\restriction B_{s(\beta)}\in x_i^\beta$ and $y\restriction B_{s(\beta)}\in x_j^\beta$. Since $\beta\in s^{-1}[\alpha+1]\cap A$, $s(\beta)\leq \alpha$. Thus $x\restriction B_{s(\beta)}=y\restriction B_{s(\beta)}$, $x_i^\beta=x_j^\beta$, and $i=j$. 
    So $F(x) (\beta)=F(y) (\beta)$ for all $\beta\in s^{-1}[\alpha+1]\cap A$. We conclude that $F(x)\restriction s^{-1}[\alpha+1] =F(y)\restriction $, hence $H(x\restriction B_\alpha)=H(y\restriction B_\alpha)$ and $H$ is well defined.

   Finally, $F$ is $A$-recursive by definition of \(H\). 
    \end{proof}
    
\ref{approx2}$\Rightarrow$ \ref{approx1}.
Let $\seq{B_\alpha}{\alpha<\kappa}$ be an approximation sequence of $\mathcal{I}$ and $f:\pre{\k}{\kappa}\rightarrow\pre{\k}{\kappa}$ be an $A$-recursive function over $\seq{B_\alpha}{\alpha<\kappa}$ witnessing $E\ \hookrightarrow_{R(A)}\ \sim^\kappa_\mathcal{F}$. 
    For every $\alpha<\kappa$, define $E_\alpha\subseteq \pre{B_{s(\alpha)}}{\k}\times\pre{B_{s(\alpha)}}{\k}$ as follows:
    \begin{itemize}
        \item if $\alpha\notin A$, let $E_\alpha=\pre{B_{s(\alpha)}}{\k}\times\pre{B_{s(\alpha)}}{\k}$ (notice that $\alpha\notin A$ implies $f(x)(\alpha)=0$);
        \item if $\alpha\in A$, let \(E_\a\) be such that $x\ E_\alpha\ y$ if and only if $f(x^0)(\alpha)=f(y^0)(\alpha)$.
   \end{itemize}
     We now show that $\langle E_\alpha\mid \alpha<\kappa\rangle$ is a precise $A$-approximation over $\seq{B_\alpha}{\alpha<\kappa}$.
It is clear that $E_\alpha$ is an equivalence relation for every $\alpha<\k$, with at most $\k$ many classes. To see that $\langle E_\alpha\mid \alpha<\kappa\rangle$ satisfies Definition~\ref{A-approximationk}\ref{A-approximationk_2}, let $x,y\in\pre{\k}{\kappa}$  be such that $x Ey$. Since $f$ is an $A$-reduction from $E$ to $\sim^\kappa_\mathcal{F}$, $f(x)\sim^\kappa_\mathcal{F} f(y)$. Thus there is $W\in \mathcal{F}$ such that $W\subseteq \{\alpha<\kappa\mid f(x)(\alpha)= f(y)(\alpha)\}$. 
Since $\mathcal{F}=\mathcal{H}\restriction A$, there is $D\in \mathcal{B}$ such that $D\cap A\subseteq\{\alpha<\kappa\mid f(x)(\alpha)= f(y)(\alpha)\}$.
On the other hand $f$ is $A$-recursive over $\seq{B_\alpha}{\alpha<\kappa}$, there is $H:\Fn_{\I}\rightarrow \Fn_{\I}$ such that for all $\alpha\in \kappa$, $$f((x\restriction B_{s(\alpha)})^0)(\alpha)=H((x\restriction B_{s(\alpha)})^0\restriction B_{s(\alpha)})(\alpha)=H(x\restriction B_{s(\alpha)})(\alpha)=f(x)(\alpha)$$ and $$f(y)(\alpha)=f((y\restriction B_{s(\alpha)})^0)(\alpha).$$ 
Since for all $\alpha\in D\cap A$, $f(x)(\alpha)=f(y)(\alpha)$, $f((x\restriction B_{s(\alpha)})^0)(\alpha)=f((y\restriction B_{s(\alpha)})^0)(\alpha)$ holds for all $\alpha\in D\cap A$. We conclude that for all $\alpha\in D\cap A$, $x\restriction B_{s(\alpha)}\ E_\alpha\ y\restriction B_{s(\alpha)}$. 

We conclude by showing that $\langle E_\alpha\mid \alpha<\kappa\rangle$ satisfies Definition~\ref{A-approximationk}\ref{A-approximationk_3}. Let $x,y\in\pre{\k}{\kappa}$ be such that $\neg(x\ E\ y)$. There is $Y\in \mathcal{F}^+$ such that for all $\alpha\in Y$, $f(x)(\alpha)\neq f(y)(\alpha)$. Since $\mathcal{F}=\mathcal{H}\restriction A$, $A\cap Y\in \mathcal{F}^+$. Let us denote $A\cap Y$ by $A'$. Therefore, $A'\subseteq A$ is such that for all $\alpha\in A'$, $f(x)(\alpha)\neq f(y)(\alpha)$. Using a similar argument as above, we conclude that $f((x\restriction B_{s(\alpha)})^0)(\alpha)\neq f((y\restriction B_{s(\alpha)})^0)(\alpha)$ holds for all $\alpha\in A'$. We conclude that for all $\alpha\in A'$, $\neg (x\restriction B_{s(\alpha)}\ E_\alpha\ y\restriction B_{s(\alpha)})$.
\end{proof}

\begin{corollary}
Let $\mathcal{H}=\{\k\}$, $E$ an equivalence relation on $\pre{\k}{\kappa}$ and let $X\subseteq \k$ be an unbounded set. Then, $E$ has a precise $X$-approximation if and only if $E\hookrightarrow_{R(X)} X_\k$. 
In particular, $E$ has a precise $\k$-approximation if and only if $E\hookrightarrow_{R(\k)} \id$.
\end{corollary}

\begin{corollary}
Let $\mathcal{H}=\CUB$, $E$ an equivalence relation on $\pre{\k}{\kappa}$, and  let $S\subseteq \k$ be a stationary set. Then, $E$ has a precise $S$-approximation if and only if $E\hookrightarrow_{R(S)}=_S^\k$. 
In particular, $E$ has a precise $\k$-approximation if and only if $E\hookrightarrow_{R(\k)} =^\k_{\CUB}$.
\end{corollary}

\section*{Acknowledgements}
The first author was supported by European Research Council (ERC) under the European Union's Horizon 2020 research and innovation programme (grant agreement No 101020762), the Research Council of Finland (decision number 368671), and the project \textit{Perspectives on computational logic}, funded by the Research Council of Finland, project number 36942.
We thank Tapani Hyttinen for his advices on analytic sets. We thank Peter Holy for the discussion on the nonstationary topology.







\end{document}